# COMPARISON OF SEMIMARTINGALES AND LÉVY PROCESSES

By Jan Bergenthum and Ludger Rüschendorf

*University of Freiburg*

In this paper, we derive comparison results for terminal values of $d$-dimensional special semimartingales and also for finite-dimensional distributions of multivariate Lévy processes. The comparison is with respect to nondecreasing, (increasing) convex, (increasing) directionally convex and (increasing) supermodular functions. We use three different approaches. In the first approach, we give sufficient conditions on the local predictable characteristics that imply ordering of terminal values of semimartingales. This generalizes some recent convex comparison results of exponential models in [*Math. Finance* **8** (1998) 93–126, *Finance Stoch.* **4** (2000) 209–222, *Proc. Steklov Inst. Math.* **237** (2002) 73–113, *Finance Stoch.* **10** (2006) 222–249]. In the second part, we give comparison results for finite-dimensional distributions of Lévy processes with infinite Lévy measure. In the first step, we derive a comparison result for Markov processes based on a monotone separating transition kernel. By a coupling argument, we get an application to the comparison of compound Poisson processes. These comparisons are then extended by an approximation argument to the ordering of Lévy processes with infinite Lévy measure. The third approach is based on mixing representations which are known for several relevant distribution classes. We discuss this approach in detail for the comparison of generalized hyperbolic distributions and for normal inverse Gaussian processes.

**1. Introduction.** Stochastic ordering and comparison results for stochastic models are topics of particular interest which have undergone an intensive development. The main aim is to derive approximation results and bounds for complex stochastic models where fundamental model assumptions are only partially known. Important fields of application include insurance mathematics, risk management, reliability theory, renewal theory, stochastic networks, financial mathematics, statistical physics and many others. Several comprehensive books survey wide areas of this research (see [15, 17, 22, 23, 24]). Most of the development is concentrated on one- or multi-dimensional orders. There is also an extended theory for ordering of discrete-time processes (like queuing sequences, renewal sequences, Markov chains) or related point processes with a wide variety of applications. Stochastic ordering results for the classical stochastic order (w.r.t. nondecreasing functions) have also been established under various conditions for diffusion-type processes (cf. [6, 7, 11, 18, 25]) and for Markov processes (cf. [16]). These results are parallel to classical comparison theorems for solutions of differential equations. Some interesting convex comparison results for exponential stochastic models have been developed in recent financial mathematics papers (see [2, 3, 5, 8, 9, 10]). The main aim in these papers is to derive sharp upper or lower bounds for option prices. The methods used in these papers are based on stochastic calculus (Itô's formula) and the propagation of convexity property (see [2, 3, 5, 8]) as well as on the coupling method (see [9, 10]).

In our paper, we first give an extension of the approach in [2, 3, 5, 8] to obtain comparison results for multidimensional semimartingales for various diffusion and monotonicity stochastic orders. This leads to stochastic ordering results of the following type. Suitable ordering of the local characteristics plus a "propagation of ordering" property imply that the basic semimartingale $S$ can be compared to a Markovian semimartingale $S^*$ at any fixed time $T$. This result is worked out in Section 2. The main tool for the proof of the comparison result is the supermartingale property of the comparison process $\mathcal{G}(t, S_t)$ based on the propagation operator $\mathcal{G}(t, s)$, which for an ordering function $g$, is given by

$$\mathcal{G}(t,s) = E^*(g(S_T^*)|S_t^* = s).$$

Our comparison result has several applications. We state several comparison results between diffusion processes, stochastic volatility models, diffusions with jumps and Lévy processes. As a very specific application, we obtain a complete list of the ordering properties of multivariate normal distributions w.r.t. diffusion and monotonicity orders.

In Section 3, we use a different technique to obtain ordering results for the finite-dimensional distributions of Lévy processes. For the comparison, we use an approximation result to reduce the ordering problem to the comparison of compound Poisson processes. These can be ordered by a direct coupling argument and by a general ordering result for time-homogeneous Markov processes based on the existence of a monotone separating transition kernel. The ordering conditions are then formulated in terms of the Lévy measures as well as the drift and diffusion coefficients. As an example, we discuss ordering properties of normal inverse Gaussian (NIG) processes.

Finally, in Section 4, we derive ordering results based on mixing-type representations which are known for several interesting classes of distributions like generalized hyperbolic (GH), multivariate $t$ and elliptically contoured distributions (cf. [4] for applications to risk management). We discuss this method in detail for the particular case of GH-distributions. Our main results describe the ordering properties of GH-distributions in terms of the parameters of the mixing variables, which are identical with the basic parameters of the distributions. For NIG-processes, these results also imply ordering results for the finite-dimensional distributions.

**2. Comparison results for semimartingales.** We derive comparisons of terminal values of multivariate special semimartingales with absolutely continuous characteristics and also of finite-dimensional distributions of Lévy processes. The comparison is based on stochastic orderings $\leq_{\mathcal{F}}$ induced by the following function classes $\mathcal{F}$:

$$
\begin{aligned}
\mathcal{F}_{\mathrm{st}} &:= \{f : \mathbb{R}^d \to \mathbb{R}, f \text{ is increasing}\}, \\
\mathcal{F}_{\mathrm{cx}} &:= \{f : \mathbb{R}^d \to \mathbb{R}, f \text{ is convex}\}, \\
\mathcal{F}_{\mathrm{dcx}} &:= \{f : \mathbb{R}^d \to \mathbb{R}, f \text{ is directionally convex}\}, \\
\mathcal{F}_{\mathrm{sm}} &:= \{f : \mathbb{R}^d \to \mathbb{R}, f \text{ is supermodular}\}, \\
\mathcal{F}_{\mathrm{icx}} &:= \mathcal{F}_{\mathrm{cx}} \cap \mathcal{F}_{\mathrm{st}}, \qquad \mathcal{F}_{\mathrm{idcx}} := \mathcal{F}_{\mathrm{dcx}} \cap \mathcal{F}_{\mathrm{st}}, \qquad \mathcal{F}_{\mathrm{ism}} := \mathcal{F}_{\mathrm{sm}} \cap \mathcal{F}_{\mathrm{st}}
\end{aligned}
$$

(1)

or by subclasses of them, in particular by single elements $g$ of these classes. We consider the componentwise order on $\mathbb{R}^d$, that is, for $x, y \in \mathbb{R}^d$, $x \leq y$ holds if $x^i \leq y^i$ for all $i \leq d$. For random vectors $X$ and $Y$, the stochastic ordering $X \leq_{\mathcal{F}} Y$ is defined by $Ef(X) \leq Ef(Y)$ for all $f \in \mathcal{F}$ such that $f(X)$ and $f(Y)$ are integrable. For the specific classes $\mathcal{F}$ in (1), we denote the orderings $X \leq_{\mathcal{F}} Y$ by $X \leq_{\mathrm{st}} Y$, $X \leq_{\mathrm{cx}} Y$, $X \leq_{\mathrm{dcx}} Y, \ldots$. Orderings induced by these function classes are also generated by all functions from $\mathcal{F} \cap C^\infty$, where $\mathcal{F}$ denotes one of the classes in (1) and $C^\infty$ is the set of infinitely differentiable functions, as well as by many other order generating function classes $\mathcal{F}_0 \subset \mathcal{F}$. For results on stochastic orders, we refer to [15, 17, 22] and [24].

For a finite time horizon $T < \infty$ and $d \in \mathbb{N}$, we consider $d$-dimensional special semimartingales $S$ and $S^*$ on stochastic bases $(\Omega, \mathcal{A}, (\mathcal{A}_t)_{t \in [0,T]}, P)$ and $(\Omega^*, \mathcal{A}^*, (\mathcal{A}_t^*)_{t \in [0,T]}, P^*)$, respectively, with predictable characteristics $(B(h), C, \nu)$ and $(B^*(h^*), C^*, \nu^*)$, respectively, where the drift components $B(h), B^*(h^*)$ depend on the choice of truncation functions $h, h^*$. We assume throughout the paper that the characteristics are absolutely continuous and denote the corresponding differential characteristics by $b(h), c$ and $K$, where $b(h)$ is a predictable $d$-dimensional process, $c = (c^{ij})_{i,j \leq d}$ is predictable with values in the set of symmetric, positive semidefinite $d \times d$ matrices with real entries $M_+(d, \mathbb{R})$ and $K_{\omega,t}(dx)$ is a transition kernel from $(\Omega \times [0,T], \mathcal{P})$ into $(\mathbb{R}^d, \mathcal{B}^d)$ that satisfies $K_{\omega,t}(0) = 0$ and $\int K_{\omega,t}(dx)(|x|^2 \wedge 1) \leq 1$ (cf. [13],

Proposition II.2.9). For $A, B \in M_+(d, \mathbb{R})$, the positive semidefinite order $A \leq_{\text{psd}} B$ is defined by $x^T(B-A)x \geq 0$ for all $x \in \mathbb{R}^d$.

By
$$S \sim (b(h), c, K) \sim (b(h), c, K)_h$$
we indicate that $S$ has predictable differential characteristics $b(h), c$ and $K$. Special semimartingales are characterized by $(|x|^2 \wedge |x|) * \nu \in \mathcal{A}_{\text{loc}}$ and $h = \text{id}$ may serve as "truncation function" (cf. [13], Proposition II.2.29). We define $b := b(\text{id})$. For example, Lévy processes with existing first moments are special semimartingales with deterministic differential characteristics that are given by

(2) $$b_t(\omega) = tb, \qquad c_t(\omega) = t\Sigma, \qquad K_t(\omega) = tF,$$

where $b \in \mathbb{R}^d$ and $\Sigma \in M(d, \mathbb{R})$ are constant and $F$ denotes the Lévy measure (cf. [13], Corollary II.4.19, [21], Corollary 5.25.8). In this case, we write $S \sim (b, \Sigma, F)$. Another example for a special semimartingale is a diffusion with jumps that is given by an SDE

$$dS_t = b(t, S_t)\,dt + \sigma(t, S_t)\,dW_t + \phi(t, S_{t-}, y)(N(dt, dy) - q(dt, dy)),$$
$$S_0 = x,$$

where $W$ is a $d$-dimensional Brownian motion and $N$ is a Poisson random measure on $\mathbb{R}_+ \times \mathbb{R}^d$ with intensity measure $q(dt, dy) = dt \otimes \lambda(dy)$, that satisfies $(|x|^2 \wedge |x|) * \nu \in \mathcal{A}_{\text{loc}}$. In this case, $S$ has characteristics $S \sim (b, c, K)$, with

(3) $$c = \sigma\sigma^T, \qquad K_t(s, A) = \int \mathbb{1}_{A\setminus\{0\}}(\phi(t, s, y))\lambda(dy) =: \lambda^{\phi(t,s,\cdot)}(A);$$

see [13], Remark III.2.29.

For the special semimartingale $S^*$, we additionally assume a Markovian structure in the sense that its differential characteristics are deterministic functions of time $t$ and $S^*_{t-}$, that is,

$$S^* \sim (b^*(t, S^*_{t-}), c^*(t, S^*_{t-}), K^*(t, S^*_{t-})).$$

For notation and results on semimartingales, we refer to [13].

The first approach relies on an integro-partial-differential equation for the propagation operator $\mathcal{G}(t, s)$. If $\mathcal{G}(t, S^*)$ is a local $(\mathcal{A}^*_t)$-martingale under $P^*$ for some $\mathcal{G} \in C^{1,2}([0, T] \times \mathbb{R}^d)$, then under an additional condition, it will be established in the following lemma that the propagation operator $\mathcal{G}(t, s)$ satisfies the integro-partial-differential equation

(4)
$$D_t\mathcal{G}(t,s) + \sum_{i \leq d} D_i\mathcal{G}(t,s)b^{*i}(t,s)$$
$$+ \tfrac{1}{2}\sum_{i,j \leq d} D^2_{ij}\mathcal{G}(t,s)c^{*ij}(t,s) + \int_{\mathbb{R}^d}(\Lambda\mathcal{G})(t,s,y)K^*_t(s,dy) = 0,$$

where $(\Lambda \mathcal{G})(t,s,y) := \mathcal{G}(t, s+y) - \mathcal{G}(t,s) - \sum_{i \le d} D_i \mathcal{G}(t,s) y^i$, $y \in \mathbb{R}^d$ and $y^i$ denotes the $i$th component. For the jump measure $\mu$ of $S$ and

(5) $$W(\omega, t, y) := \Lambda \mathcal{G}(t, S_{t-}(\omega), y),$$

the integral process $(W * \mu)_t$ is defined by

$$(W * \mu)_t(\omega) := \int_{[0,t] \times \mathbb{R}^d} W(\omega, u, y) \mu(\omega; du, dy)$$
$$= \int_{[0,t] \times \mathbb{R}^d} \Lambda \mathcal{G}(u, S_{u-}(\omega), y) \mu(\omega; du, dy).$$

The predictable function $W^*$ and the integral process $(W^* * \mu^*)_t$ are defined similarly. The proof of the following lemma is given in the Appendix:

LEMMA 2.1 (Kolmogorov backward equation). *Let $S_t^* \sim (b^*(t, S_{t-}^*), c^*(t, S_{t-}^*), K^*(t, S_{t-}^*))$ be a $d$-dimensional special semimartingale, let $\mathcal{G} \in C^{1,2}([0,T] \times \mathbb{R}^d)$ and assume that $\mathcal{G}(t, S_t^*)$ is a local $(\mathcal{A}_t^*)$-martingale under $P^*$. If $(|W^*| * \mu^*)_t \in \mathcal{A}_{\text{loc}}^+$ or if $\mathcal{G}(t, \cdot)$ is convex, then $\mathcal{G}(t,s)$ satisfies the Kolmogorov backward equation* (4).

In the sequel, for $g \in \mathcal{F}$, we consider the functional $\mathcal{G}(t,s)$ defined by the propagation operator

(6) $$\mathcal{G}(t,s) = E^*(g(S_T^*) | S_t^* = s),$$

where $E^*$ denotes the expectation with respect to $P^*$. A crucial assumption in the comparison results for special semimartingales is the propagation of order property: for some ordering function $g = \mathcal{G}(T, \cdot) \in \mathcal{F}$ of interest, we assume that the ordering at time $T$ is propagated to earlier times $t \in [0, T]$, that is, $\mathcal{G}(t, \cdot) \in \mathcal{F}$. This order-propagating property is known in exponential diffusion models in financial context for $\mathcal{F} = \mathcal{F}_{\text{cx}}$ as "propagation of convexity." The no-arbitrage price process $\mathcal{G}(t,s)$ of a European contingent claim with convex payoff function $g(s) = \mathcal{G}(T, s)$ is convex in $s$, for all $t \in [0, T]$ (cf. [2, 3, 5, 8, 10]).

ASSUMPTION PO($g$) (*Propagation of order*). For some $g \in \mathcal{F}$, we say that $S^*$ satisfies the propagation of order property PO($g$) if $\mathcal{G}(t, \cdot) \in \mathcal{F}$ for $0 \le t \le T$. Similarly, $S^*$ satisfies PO($\mathcal{F}_0$) for some $\mathcal{F}_0 \subset \mathcal{F}$ if PO($g$) holds for all $g \in \mathcal{F}_0$.

In the following, we need some integrability properties of $\mathcal{G}(t, S_t)$. We define the classes $\overline{\mathcal{F}} = \overline{\mathcal{F}}(S, S^*)$ and $\widetilde{\mathcal{F}} = \widetilde{F}(S, S^*)$ by

$$\overline{\mathcal{F}} := \{g \in \mathcal{F} : \mathcal{G}(t,s) \in C^{1,2}([0,T] \times \mathbb{R}^d),$$

$$\mathcal{G}(t, S_t) \text{ bounded below and } E\mathcal{G}(t, S_t) < \infty \ \forall t \in [0, T]\},$$

$$\widetilde{\mathcal{F}} := \{g \in \mathcal{F} : \mathcal{G}(t, s) \in C^{1,2}([0, T] \times \mathbb{R}^d) \text{ and}$$

$$\mathcal{G}(t, S_t) \text{ is a process of class } (D)\}.$$

If, for example, $g$ has bounded linear growth, $|g(x)| \leq K|x|$ for $K < \infty$ and $S$, $S^*$ are Lévy processes or exponential Lévy processes, then the integrability condition $E\mathcal{G}(t, S_t) = \int \int g(s + s^*) P^{*S^*_{T-t}}(ds^*) P^{S_t}(ds) < \infty$ follows from the integrability of $S_t$ and $S^*_{T-t}$.

The general comparison result is formulated in the following theorem for the increasing directionally convex and for the increasing convex orders generated by $\mathcal{F} = \overline{\mathcal{F}}_{\mathrm{idcx}}$ and $\mathcal{F} = \overline{\mathcal{F}}_{\mathrm{icx}}$. For the other order-generating function classes in (1), similar results hold true (see Remark 2.3). The proof of the following theorem is based on ideas similar to those in [2, 3, 5] and [8], which are concerned with stochastic exponential models.

THEOREM 2.2 [Increasing (directionally) convex comparison of semimartingales]. *Let $S, S^* \in \mathcal{S}_p^d$ and $S_0 = S_0^*$ have differential local characteristics $S \sim (b, c, K)_{\mathrm{id}}$ and $S^* \sim (b^*, c^*, K^*)_{\mathrm{id}}$. Let $W$ be as in (5) and let $W^*$ be defined similarly.*

1. *Let $g \in \overline{\mathcal{F}}_{\mathrm{idcx}}$ or $g \in \widetilde{\mathcal{F}}_{\mathrm{idcx}}$, let Assumption $\mathrm{PO}(g)$ hold and assume that the associated processes $(|W| * \mu)_t$ and $(|W^*| * \mu^*)_t$ are in $\mathcal{A}^+_{\mathrm{loc}}$. Then the comparison of the differential characteristics*

(7) $$b_t^i(\omega) \leq b^{*i}(t, S_{t-}(\omega)), \qquad c_t^{ij}(\omega) \leq c^{*ij}(t, S_{t-}(\omega)),$$

(8) $$\int_{\mathbb{R}^d} f(t, S_{t-}(\omega), x) K_{\omega, t}(dx) \leq \int_{\mathbb{R}^d} f(t, S_{t-}(\omega), x) K_t^*(S_{t-}(\omega), dx),$$

$\lambda \times P$-*a.e. for all $f \in \mathbb{R}_+ \times \mathbb{R}^d \times \mathbb{R}^d \to \mathbb{R}$ with $f(t, s, \cdot) \in \mathcal{F}_{\mathrm{dcx}}$ such that the integrals exist, implies that*

$$Eg(S_T) \leq E^* g(S_T^*).$$

2. *Let $g \in \overline{\mathcal{F}}_{\mathrm{icx}}$ and let Assumption $\mathrm{PO}(g)$ hold true. Then the comparison of the differential characteristics*

(9) $$b_t^i(\omega) \leq b^{*i}(t, S_{t-}(\omega)), \qquad (c_t^{ij}(\omega))_{i,j \leq d} \leq_{\mathrm{psd}} (c^{*ij}(t, S_{t-}(\omega)))_{i,j \leq d},$$

(10) $$\int_{\mathbb{R}^d} f(t, S_{t-}(\omega), x) K_{\omega, t}(dx) \leq \int_{\mathbb{R}^d} f(t, S_{t-}(\omega), x) K_t^*(S_{t-}(\omega), dx),$$

$\lambda \times P$-*a.e., for all $f \in \mathbb{R}_+ \times \mathbb{R}^d \times \mathbb{R}^d \to \mathbb{R}$ with $f(t, s, \cdot) \in \mathcal{F}_{\mathrm{cx}}$, such that the integrals exist, implies that*

$$Eg(S_T) \leq E^* g(S_T^*).$$

PROOF. We prove that the comparison process $\mathcal{G}(t, S_t)$ is an $(\mathcal{A}_t)$-supermartingale under $P$. Then using the fact that $S_0 = S_0^*$, we obtain the stated comparison result

$$Eg(S_T) = E\mathcal{G}(T, S_T) \leq \mathcal{G}(0, S_0) = E^* g(S_T^*).$$

1. Let $g \in \overline{\mathcal{F}}_{\text{idcx}}$ or $\widetilde{\mathcal{F}}_{\text{idcx}}$. As $(|W| * \mu)_t \in \mathcal{A}_{\text{loc}}^+$, Itô's formula implies, similarly as in the proof of Lemma 2.1, that

$$\begin{aligned}
\mathcal{G}(t, S_t) = {} & \mathcal{G}(0, S_0) + M_t \\
& + \int_{[0,t]} \bigg\{ D_t \mathcal{G}(u, S_{u-}) + \sum_{i \leq d} D_i \mathcal{G}(u, S_{u-}) b_u^i \\
& \qquad + \tfrac{1}{2} \sum_{i,j \leq d} D_{ij}^2 \mathcal{G}(u, S_{u-}) c_u^{ij} + \int_{\mathbb{R}^d} (\Lambda \mathcal{G})(u, S_{u-}, x) K_u(dx) \bigg\} du,
\end{aligned}$$

where $M_t := \sum_{i \leq d} \int_{[0,t]} D_i \mathcal{G}(u, S_{u-}) \, dN_u^i$ is a local $(\mathcal{A}_t)$-martingale under $P$ and $N_t^i$ denotes the $i$th component of the martingale part of the canonical decomposition of $S_t = S_0 + N_t + B_t$. As $\mathcal{G}(t, s)$ satisfies the integro-partial-differential equation (4), we obtain

$$\begin{aligned}
\mathcal{G}(t, S_t) = {} & \mathcal{G}(0, S_0) + M_t \\
& + \int_{[0,t]} \bigg\{ \sum_{i \leq d} D_i \mathcal{G}(u, S_{u-})(b_u^i - b^{*i}(u, S_{u-})) \\
& \qquad + \tfrac{1}{2} \sum_{i,j \leq d} D_{ij}^2 \mathcal{G}(u, S_{u-})(c_u^{ij} - c^{*ij}(u, S_{u-})) \\
& \qquad + \int_{\mathbb{R}^d} (\Lambda \mathcal{G})(u, S_{u-}, x)(K_u(dx) - K_u^*(S_{u-}, dx)) \bigg\} du \\
=: {} & \mathcal{G}(0, S_0) + M_t + V_t,
\end{aligned} \quad (11)$$

where $V_t$ is a predictable process of finite variation and, therefore, $V_t \in \mathcal{A}_{\text{loc}}$. We now prove that inequalities (7) and (8) imply that $\mathcal{G}(t, S_t)$ is a local $(\mathcal{A}_t)$-supermartingale under $P$. By Assumption PO($g$), $\mathcal{G}(t, \cdot)$ is increasing and directionally convex for all $t \in [0, T]$. Therefore, $D_i \mathcal{G}(t, s)$ and $D_{ij}^2 \mathcal{G}(t, s)$ are nonnegative and inequalities (7) imply that the first and second terms of the integrand of $V_t$ are nonpositive. For fixed $(\omega, u) \in (\Omega, [0, T])$, we define the function

$$\Upsilon(x) := (\Lambda \mathcal{G})(u, S_{u-}(\omega), x). \quad (12)$$

As $D_{ij}^2 \Upsilon(x) = D_{ij}^2 \mathcal{G}(u, S_{u-} + x) \geq 0$, for all $x$, $\Upsilon$ is directionally convex and, therefore, (8) implies that the last term of the integrand of $V_t$ is nonpositive. This yields $-V_t \in \mathcal{A}_{\text{loc}}^+$ and it follows that $\mathcal{G}(t, S_t)$ is a local $(\mathcal{A}_t)$-supermartingale under $P$.

It remains to prove that $\mathcal{G}(t, S_t)$ is an $(\mathcal{A}_t)$-supermartingale. If $g \in \overline{\mathcal{F}}_{\text{idcx}}$ is bounded below, $M_t$ is bounded below and, therefore, is an $(\mathcal{A}_t)$-supermartingale under $P$. It follows that $\mathcal{G}(t, S_t)$ is a supermartingale as it is integrable. In the case where $g \in \widetilde{\mathcal{F}}_{\text{idcx}}$, we consider a localizing sequence $\tau_n$ for $\mathcal{G}(t, S_t)$. Since for all $t \in [0, T]$, we have $P$-a.s. that $(\mathcal{G}(t, S_t))^{\tau_n} \to \mathcal{G}(t, S_t)$, $n \to \infty$, and $\mathcal{G}(t, S_t)$ is of class $(D)$, the convergence takes place in $L^1$ and, therefore, $\mathcal{G}(t, S_t)$ is an $(\mathcal{A}_t)$-supermartingale under $P$.

2. Let $g \in \overline{\mathcal{F}}_{\text{icx}}$. Similarly to the proof of the first part, we show that $\mathcal{G}(t, S_t)$ is an $(\mathcal{A}_t)$-supermartingale under $P$. The evolution of $\mathcal{G}(t, S_t)$ under $P$ is given in (11). Observe that as in Lemma 2.1, convexity of $\mathcal{G}(t, \cdot)$ implies $(|W| * \mu)_t, (|W^*| * \mu^*)_t \in \mathcal{A}_{\text{loc}}^+$. As $\mathcal{G}(t, \cdot)$ is increasing, comparison of the drift characteristics in (9) implies that the first term of the integrand of $V_t$ is nonpositive. Due to the positive semidefiniteness of $(c^{*ij} - c^{ij})_{i,j \leq d} := (c^{*ij}(u, S_{u-})(\omega) - c_u^{ij}(\omega))_{i,j \leq d}$ for fixed $(\omega, u)$, its spectral decomposition is given by $(\sum_{k \leq d} \lambda_k e_k^i e_k^j)_{i,j \leq d}$, where the eigenvalues $\lambda_k$ are nonnegative and $e_k = (e_k^1, \ldots, e_k^d)$ denote the eigenvectors. Therefore, the second term of the integrand of $V_t$ is of the form $-\frac{1}{2} \sum_{k \leq d} \lambda_k \sum_{i,j \leq d} D_{ij}^2 \mathcal{G}(u, S_{u-}) e_k^i e_k^j$ and is nonpositive due to the convexity of $\mathcal{G}(t, \cdot)$ from Assumption PO($g$). Also, by convexity of $\mathcal{G}(t, \cdot)$, it follows that $\Upsilon(x)$ from (12) is nonnegative and convex. Assumption (10) implies that the last term of the integrand of $V_t$ is nonpositive and, therefore, $-V_t \in \mathcal{A}_{\text{loc}}^+$. This implies that $\mathcal{G}(t, S_t)$ is a *local* $(\mathcal{A}_t)$-supermartingale under $P$. As $\mathcal{G}(t, S_t)$ is bounded below and $E\mathcal{G}(t, S_t) < \infty$, we obtain, as in the proof of the first part, that $\mathcal{G}(t, S_t)$ is an $(\mathcal{A}_t)$-supermartingale under $P$. □

REMARK 2.3. 1. As seen in the proof, it suffices to assume that

$$\int_{\mathbb{R}^d} (\Lambda \mathcal{G})(t, S_{t-}, x) K_t(dx) \leq \int_{\mathbb{R}^d} (\Lambda \mathcal{G})(t, S_{t-}, x) K_t^*(S_{t-}, dx)$$

instead of (8) and (10) and, further, that $S_0 \leq S_0^*$.

2. Similar comparison results are obtained for $g \in \widetilde{\mathcal{F}}$ (or $g \in \overline{\mathcal{F}}$), where $\mathcal{F} \in \{\mathcal{F}_{\text{st}}, \mathcal{F}_{\text{dcx}}, \mathcal{F}_{\text{cx}}, \mathcal{F}_{\text{ism}}, \mathcal{F}_{\text{sm}}\}$. The crucial equation in the proof is (11). Table 1 lists sufficient conditions similar to (7)–(10) that imply $Eg(S_T) \leq E^* g(S_T^*)$ under the respective integrability conditions.

3. For $g \in \widetilde{\mathcal{F}}$, $\mathcal{F} \in \{\mathcal{F}_{\text{st}}, \mathcal{F}_{\text{idcx}}, \mathcal{F}_{\text{ism}}, \mathcal{F}_{\text{icx}}, \mathcal{F}_{\text{dcx}}, \mathcal{F}_{\text{sm}}, \mathcal{F}_{\text{cx}}\}$, similar arguments show that $E^* g(S_T^*) \leq Eg(S_T)$, if the corresponding inequalities of the differential characteristics of $S$ and $S^*$ are reversed.

The propagation of order property is a crucial assumption in the previous comparison theorem. To apply the theorem, this condition needs to be checked in the models considered. Bergenthum and Rüschendorf [3] established propagation of (directional) convexity for stochastic exponentials of Lévy processes and the propagation of convexity property for some classes of diffusions and diffusions with jumps with zero drift coefficient and finite jump intensity. We establish propagation of increasing convexity for some classes of diffusions with jumps with nonzero drift and finite jump intensity and in the second part of this section, we establish the PO($\mathcal{F}$) property for spatially homogeneous processes for all order-generating function classes $\mathcal{F}$ considered in this paper.

Let $S^*$ be a unique strong solution of the SDE

$$\begin{aligned}dS_t^* &= b^*(t, S_t^*)\, dt + \sigma^*(t, S_t^*)\, dW_t^* \\ &\quad + \phi^*(t, S_{t-}^*, y)(N^*(dt, dy) - \lambda^*(dy)\, dt),\\ S_0^* &= s_0,\end{aligned} \tag{13}$$

where $W_t^*$ is a $d$-dimensional Brownian motion and $N^*$ is a Poisson random measure on $[0, T] \times \mathbb{R}$ with deterministic intensity $\lambda^*(dy)\, dt$ and $\lambda^*(\mathbb{R}) < \infty$. We consider the following two cases for the parameter of the jump part $\phi^*$ of the SDE (13). In the first case, we assume that $\phi^*$ factorizes as

$$\phi^*(t, s, y) := \varphi^*(t, s) \psi^*(t, y), \tag{14}$$

where $\varphi^* : [0, T] \times \mathbb{R}^d \to \mathbb{R}_+$, $\varphi^*(t, \cdot)$ is increasing and convex for all $t \in [0, T]$ and $\psi^* : [0, T] \times \mathbb{R} \to \mathbb{R}^d$. In the second case, we assume that $d = 1$ and $\phi^*$ is of the form

$$\phi^*(t, s, y) := \sum_{i \leq m} \varphi_i^*(t, s) \psi_i^*(t, y), \tag{15}$$

where $\varphi_i^* : [0, T] \times \mathbb{R} \to \mathbb{R}_+$, $\varphi_i^*(t, \cdot) \in \mathcal{F}_{\mathrm{icx}}$ for all $t \in [0, T]$ and $\psi_i^* : [0, T] \times \mathbb{R} \to \mathbb{R}_+$, $\psi_i^*(t, \cdot)$ is increasing for all $t \in [0, T]$, for all $i \leq m$.

TABLE 1

| Ordering | Drift | Diffusion | Jump |
| --- | --- | --- | --- |
| $\mathcal{F}_{\mathrm{st}}$ | $b \leq b^*$ | $c = c^*$ | $K = K^*$ |
| $\mathcal{F}_{\mathrm{idcx}}$ | $b \leq b^*$ | $c \leq c^*$ | $K \leq_{\mathrm{dcx}} K^*$ |
| $\mathcal{F}_{\mathrm{icx}}$ | $b \leq b^*$ | $c \leq_{\mathrm{psd}} c^*$ | $K \leq_{\mathrm{cx}} K^*$ |
| $\mathcal{F}_{\mathrm{ism}}$ | $b \leq b^*$ | $c \leq c^*, c^{ii} = c^{*ii}, i \leq d$ | $K \leq_{\mathrm{sm}} K^*$ |
| $\mathcal{F}_{\mathrm{dcx}}$ | $b = b^*$ | $c \leq c^*$ | $K \leq_{\mathrm{dcx}} K^*$ |
| $\mathcal{F}_{\mathrm{cx}}$ | $b = b^*$ | $c \leq_{\mathrm{psd}} c^*$ | $K \leq_{\mathrm{cx}} K^*$ |
| $\mathcal{F}_{\mathrm{sm}}$ | $b = b^*$ | $c \leq c^*, c^{ii} = c^{*ii}, i \leq d$ | $K \leq_{\mathrm{sm}} K^*$ |

REMARK 2.4. Under the assumption in (14), the SDE in (13) can be seen as a diffusion driven by a process with independent increments (PII). For time-independent $\psi^*$, it is a Lévy-driven diffusion. In the case $d = 1$, the class of functions in (15) allows the approximation of general linking functions $\phi^*$ and, thus, in this case, the assumption allows the treatment of more general diffusions with jumps.

For $t_0 \in [0, T]$ and $K \in \mathbb{N}$, we discretize $[t_0, T]$ into $K + 1$ equidistant points $t_i := i\frac{T-t_0}{K} + t_0, i \in \{0, \ldots, K\}$ and denote the Euler scheme $\widetilde{S}_K^*$ of $S^*$ by

$$\widetilde{S}_{K,t_{i+1}}^* = \widetilde{S}_{K,t_i}^* + b^*(t_i, \widetilde{S}_{K,t_i}^*)\Delta t_i + \sigma^*(t_i, \widetilde{S}_{K,t_i}^*)(W_{t_{i+1}}^* - W_{t_i}^*)$$
(16)
$$+ \phi^*(t_i, \widetilde{S}_{K,t_i}^*, Y^*)\widetilde{N}^* - E^{Y^*}\phi^*(t_i, \widetilde{S}_{K,t_i}^*, Y^*)\lambda^*(\mathbb{R})\Delta t_i,$$

$$\widetilde{S}_{K,t_0}^* = s,$$

where $\widetilde{N}^*$ is binomial with $P^{\widetilde{N}^*} = (1 - \lambda^*(\mathbb{R})\Delta t_i)\varepsilon_{\{0\}} + \lambda^*(\mathbb{R})\Delta t_i\varepsilon_{\{1\}}$, $Y^*$ has distribution $\frac{\lambda^*(dy)}{\lambda^*(\mathbb{R})}$ on $(\mathbb{R}, \mathcal{B})$, $E^{Y^*}f(\cdot, \cdot, Y^*) := \frac{1}{\lambda^*(\mathbb{R})}\int f(\cdot, \cdot, y)\lambda^*(dy)$ and $\Delta t_i = t_{i+1} - t_i = \frac{T-t_0}{K}$. Observe that $W_{t_{i+1}}^* - W_{t_i}^* \sim N(0, \Delta I)$, where $N(\mu, \Sigma)$ denotes the $d$-dimensional normal distribution with expectation $\mu$ and covariance matrix $\Sigma$ and $I \in M_+(d, \mathbb{R})$ is the identity.

We establish propagation of increasing convexity for the propagation operator $\widetilde{\mathcal{G}}_K(t, s)$ of the Euler scheme $\widetilde{S}_K^*$ that is given by

$$\widetilde{\mathcal{G}}_K(t, s) := E^*(g(\widetilde{S}_{K,T}^*)|\widetilde{S}_{K,t}^* = s).$$

We then use an approximation argument that implies propagation of increasing convexity for the propagation operator $\mathcal{G}(t, s)$ also. We call this the *approximation property*.

ASSUMPTION AP($g$) (*Approximation property*). Let $g : \mathbb{R}^d \to \mathbb{R}$. The Euler scheme $\widetilde{S}_K^*$ of $S^*$ satisfies the approximation property AP($g$), if for $K \to \infty$, it holds true that

$$\widetilde{\mathcal{G}}_K(t, s) := E^*(g(\widetilde{S}_{K,T}^*)|\widetilde{S}_{K,t}^* = s) \to \mathcal{G}(t, s) \qquad \text{for all } t \in [0, T], s \in \mathbb{R}^d.$$

Similarly, the Euler scheme $\widetilde{S}_K^*$ of $S^*$ satisfies the approximation property AP($\mathcal{F}_0$) for some $\mathcal{F}_0 \subset \mathcal{F}$ if AP($g$) holds true for all $g \in \mathcal{F}_0$.

The approximation property is satisfied under smoothness and linear growth conditions on the coefficients and smoothness and growth conditions on the test function $g$ (cf., e.g., [14]). For suitable convergence results for Lévy-driven diffusions, we refer to [12] and [19].

Sufficient for propagation of increasing convexity for the propagation operator $\widetilde{\mathcal{G}}_K(t,s)$ is $\leq_{\text{icx}}$-monotonicity of the corresponding transition operator $\mathcal{T}^*$ (i.e., $S_1^* \leq_{\text{icx}} S_2^*$ implies $\mathcal{T}^* S_1^* \leq_{\text{icx}} \mathcal{T}^* S_2^*$) that is of the form

$$\begin{aligned}\mathcal{T}^* S^* &\stackrel{d}{=} S^* + \Delta b^*(S^*) + \sigma^*(S^*) W^* \\ &\quad + \phi^*(S^*, Y^*) N^* - E^{Y^*} \phi^*(S^*, Y^*) E^* N^*.\end{aligned} \tag{17}$$

LEMMA 2.5 (Increasing convex monotonicity of Markov operators). *Let $S^*, W^*, N^*, Y^*$ be independent random variables, where $S^*$, $W^*$ are $\mathbb{R}^d$-valued, $W^* \sim N(0, \Delta I)$, $\Delta \in \mathbb{R}_+$ and $N^*, Y^*$ have values in $\mathbb{R}$. Assume that $b^*: \mathbb{R}^d \to \mathbb{R}^d$ and $\sigma^*: \mathbb{R}^d \to M_+(d, \mathbb{R})$ are increasing and convex, where $M_+(d, \mathbb{R})$ is supplied with the positive semidefinite ordering $\leq_{\text{psd}}$. Then the transition operator $\mathcal{T}^*$ given by (17) is $\leq_{\text{icx}}$-monotone if:*

1. $\phi^*(s,y) = \varphi^*(s)\psi^*(y)$, where $\varphi^*: \mathbb{R}^d \to \mathbb{R}_+$ is increasing and convex and $\psi^*: \mathbb{R} \to \mathbb{R}^d$, or
2. $d = 1$, $N^* \geq 0$ and $\phi^*(s,y) = \sum_{i \leq m} \varphi_i^*(s) \psi_i^*(y)$, where $\varphi_i^*: \mathbb{R} \to \mathbb{R}_+$ are increasing and convex and $\psi_i^*: \mathbb{R} \to \mathbb{R}_+$ are increasing for all $i \leq n$.

PROOF. Assume that $S_1^*$ and $S_2^*$ are $d$-dimensional random vectors that are independent of $W^*, N^*$ and $Y^*$ and satisfy $S_1^* \leq_{\text{icx}} S_2^*$. Without loss of generality, by Strassen's theorem, we choose $S_1^*, S_2^*$ such that $E(S_2^*|S_1^*) \geq S_1^*$. For $f \in \mathcal{F}_{\text{icx}}$, it follows from Jensen's inequality and independence of the random vectors that

$$\begin{aligned}E^* f(\mathcal{T}^* S_2^*) &= E^* E^*(f(\mathcal{T}^* S_2^*)|S_1^*, W^*, N^*, Y^*) \\ &\geq E^* f(E^*(S_2^*|S_1^*) + \Delta E^*(b^*(S_2^*)|S_1^*) \\ &\quad + E^*(\sigma^*(S_2^*)|S_1^*) W^* + E^*(\phi^*(S_2^*, Y^*)|S_1^*, Y^*) N^* \\ &\quad - E^*(E^{Y^*} \phi^*(S_2^*, Y^*)|S_1^*) E^* N^*).\end{aligned} \tag{18}$$

If $\phi^*$ factorizes into $\phi^*(s,y) = \varphi^*(s)\psi^*(y)$ as in part 1, then using the independence properties, the right-hand side of (18) equals

$$\int E^* f(E^*(S_2^*|s_1) + \Delta E^*(b^*(S_2^*)|s_1) + E^*(\sigma^*(S_2^*)|s_1) W^* \\ + E^*(\varphi^*(S_2^*)|s_1)(\psi^*(Y^*) N^* - E^*(\psi^*(Y^*) N^*))) P^{S_1^*}(ds_1).$$

Due to monotonicity and convexity of $b^*$ and $\sigma^*$, Jensen's inequality implies that $E^*(b^*(S_2^*)|s_1) \geq b^*(E^*(S_2^*|s_1)) \geq b^*(s_1)$ and $\Sigma^*(s_1) := E^*(\sigma^*(S_2^*)|s_1) \geq_{\text{psd}} \sigma^*(E^*(S_2^*|s_1)) \geq_{\text{psd}} \sigma^*(s_1)$. As $\sigma^*(s_1)$ and $\Sigma^*(s_1)$ are positive semidefinite, it holds true that $(\sigma^*(s_1))^T \sigma^*(s_1) \leq_{\text{psd}} (\Sigma^*(s_1))^T \Sigma^*(s_1)$. A classical convex ordering result for normally distributed random vectors (cf. [17], Theorem

3.4.6) implies that $\sigma^*(s_1)W^* \leq_{\text{cx}} \Sigma^*(s_1)W^*$. In the following, we establish the ordering of the jump parts

$$B := \varphi^*(s_1)(\psi^*(Y^*)N^* - E^*(\psi^*(Y^*)N^*))$$
$$\leq_{\text{cx}} E^*(\varphi^*(S_2^*)|s_1)(\psi^*(Y^*)N^* - E^*(\psi^*(Y^*)N^*)) =: C.$$

As $\varphi^* \in \mathcal{F}_{\text{icx}}$, Jensen's inequality implies that $\vartheta^*(s_1) := E^*(\varphi^*(S_2^*)|s_1) - \varphi^*(s_1) \geq 0$. For $j \leq d$, we define $R_j = \vartheta^*(s_1)(\psi_j^*(Y^*)N^* - E^*(\psi_j^*(Y^*)N^*))$. Then for $g \in \mathcal{F}_{\text{cx}} \cap C^2$, it holds true that

$$E^*g(C) \geq E^*g(B) + E^*\langle \nabla g(B), R \rangle,$$

where $R = (R_1, \ldots, R_d)$, $\nabla$ is the gradient and $\langle \cdot, \cdot \rangle$ stands for the scalar product in $\mathbb{R}^d$. From $R = \frac{\vartheta^*(s_1)}{\varphi^*(s_1)}B$ it follows that $E^*\langle \nabla g(B), R \rangle = \frac{\vartheta^*(s_1)}{\varphi^*(s_1)} E^*\langle \nabla g(B), B \rangle$. A characterization result of optimal couplings in [20], Theorem 1, implies that $(B, \nabla g(B))$ is an optimal $\ell_2$-coupling, hence $E^*\langle \nabla g(B), B \rangle \geq \langle E^*\nabla g(B), E^*B \rangle = 0$, as $E^*B = 0$. As the convex order is generated by functions $g \in \mathcal{F}_{\text{icx}} \cap C^2$, it follows that $B \leq_{\text{cx}} C$, hence $B \leq_{\text{icx}} C$.

As the increasing convex order is stable under convolutions, it follows from the independence of $W^*$, $Y^*$ and $N^*$ that

$$s_1 + \Delta b^*(s_1) + \sigma^*(s_1)W^* + \varphi^*(s_1)(\psi^*(Y^*)N^* - E^*(\psi^*(Y^*)N^*))$$
$$\leq_{\text{icx}} E^*(S_2^*|s_1) + \Delta E^*(b(S_2^*)|s_1) + E^*(\sigma(S_2^*)|s_1)W^*$$
$$+ E^*(\varphi^*(S_2^*)|s_1)(\psi^*(Y^*)N^* - E^*(\psi^*(Y^*)N^*))$$

and as $E^*f(\mathcal{T}^*S_1^*) = \int E^*f(s_1 + \Delta b^*(s_1) + \sigma^*(s_1)W^* + \varphi^*(s_1)(\psi^*(Y^*)N^* - E^*(\psi^*(Y^*)N^*)))P^{S_1^*}(ds_1)$, the result follows.

If $d = 1$ and $\phi^*$ has a representation of the form $\phi^*(s, y) = \sum_{i \leq m} \varphi_i^*(s)\psi_i^*(y)$ as in part 2, then similarly to the previous case, it suffices to establish

$$B := \sum_{i \leq m} \varphi_i^*(s_1)(\psi_i^*(Y^*)N^* - E^*(\psi_i^*(Y^*)N^*))$$
$$\leq_{\text{cx}} \sum_{i \leq m} E^*(\varphi_i^*(S_2^*)|s_1)(\psi_i^*(Y^*)N^* - E^*(\psi_i^*(Y^*)N^*)) =: C.$$

From Jensen's inequality and monotonicity of $\varphi_i^*$, it follows that $\vartheta_i^*(s_1) := E^*(\varphi_i^*(S_2^*)|s_1) - \varphi_i^*(s_1) \geq 0$. For $R := \sum_{i \leq m} \vartheta_i^*(\psi_i^*(Y^*)N^* - E^*(\psi_i^*(Y^*)N^*))$ and $g \in \mathcal{F}_{\text{cx}} \cap C^2$, Taylor's formula implies

$$E^*g(C) \geq E^*g(B) + E^*g'(B)R$$

and it remains to verify nonnegativity of $E^*g'(B)R$. We make use of some results on association of random vectors (cf. [17], Theorems 3.10.5, 3.10.7). As $Y^*$ and $N^*$ are independent random variables, $(Y^*, N^*)$ is associated.

From monotonicity of $\psi_i^* \geq 0$, it follows that $\Psi_i^*(y, n) := \psi_i^*(y)n$, $n \geq 0$, is nondecreasing in $(y, n)$, for all $i \leq m$ and

$$Z^* := (Z_1^*, \ldots, Z_m^*) = (\Psi_1^*(Y^*, N^*), \ldots, \Psi_m^*(Y^*, N^*))$$

is associated, as $N^* \geq 0$. Thus, $\overline{Z}^* = Z^* - E^*Z^*$ is associated (cf. [17], Theorem 3.10.7). From nonnegativity of $\varphi_i^*(s_1)$ and $\vartheta_i^*(s_1)$, it follows that $(B, R) = (\sum_{i \leq m} \varphi_i^*(s_1)\overline{Z}_i^*, \sum_{i \leq m} \vartheta_i^*(s_1)\overline{Z}_i^*)$ is nondecreasing in $\overline{Z}^*$ and, therefore, is associated. Thus, it holds true that $E^*F_1(B, R)F_2(B, R) \geq E^*F_1(B, R)E^*F_2(B, R)$ for all nondecreasing $F_k : \mathbb{R} \times \mathbb{R} \to \mathbb{R}$, $k = 1, 2$. As $g \in \mathcal{F}_{cx} \cap C^2$, $g'$ is nondecreasing and with $F_1(B, R) := g'(B)$, $F_2(B, R) := R$, it follows that $E^*g'(B)R \geq E^*g'(B)E^*R = 0$. $\square$

From Lemma 2.5 and Assumption AP($g$), propagation of increasing convexity for some classes of diffusions with jumps follows.

THEOREM 2.6 (Propagation of increasing convexity for diffusions with jumps). *Let $g \in \mathcal{F}_{icx}$ and $S^*$ be a $d$-dimensional diffusion with jumps. Let Assumption AP($g$) be satisfied. If $b^*(t, \cdot)$ and $\sigma^*(t, \cdot)$ are increasing and convex for all $t \in [0, T]$, as in Lemma 2.5, and $\phi^*$ satisfies (14) or (15) ($d = 1$ in the latter case), then the propagation of increasing convexity property holds, that is,*

$$\mathcal{G}(t, \cdot) \in \mathcal{F}_{icx} \qquad \forall t \in [0, T].$$

PROOF. We consider the Euler approximation scheme $\widetilde{S}_K^*$ defined in (16) with interpolation points $t_i$ and define the corresponding transition operator by $\mathcal{T}_{t_i}S^* \stackrel{d}{=} S^* + b^*(t_i, S^*)\Delta t_i + \sigma^*(t_i, S^*)W_i^* + \phi^*(t_i, S^*, Y)N_i^* - E^Y\phi^*(t_i, S^*, Y)EN_i^*$, where $W_i^* \stackrel{d}{=} W_{t_{i+1}}^* - W_{t_i}^*$ and $P^*(N_i^* = 1) = \lambda(\mathbb{R})\Delta t_i = 1 - P^*(N_i^* = 0)$. Then for $t_0 \in [0, T]$ the Markov property implies

$$\widetilde{\mathcal{G}}_K(t_0, y) = E^*(g(\widetilde{S}_{K,T}^*)|\widetilde{S}_{K,t_0}^* = y) = E^*g(\mathcal{T}_{t_{K-1}} \cdots \mathcal{T}_{t_0} y).$$

For $y_1, y_2 \in \mathbb{R}^d$ and $\alpha \in (0, 1)$, let $Y$ be a Bernoulli random vector with distribution $P^Y = \alpha\varepsilon_{\{y_1\}} + (1 - \alpha)\varepsilon_{\{y_2\}}$. Then

$$\alpha y_1 + (1 - \alpha)y_2 = EY \leq_{cx} Y.$$

Using the $\leq_{icx}$-monotonicity of the operator $\mathcal{T}_t$ for all $t \in [0, T]$, from Lemma 2.5, we obtain

$$\widetilde{\mathcal{G}}_K(t_0, \alpha y_1 + (1 - \alpha)y_2) = \widetilde{\mathcal{G}}_K(t_0, EY) = E^*g(\mathcal{T}_{t_{K-1}} \cdots \mathcal{T}_{t_0} EY)$$
$$\leq E^*g(\mathcal{T}_{t_{K-1}} \cdots \mathcal{T}_{t_0} Y) = \widetilde{\mathcal{G}}_K(t_0, Y).$$

Taking expectations on both sides implies convexity of $\widetilde{\mathcal{G}}_K(t_0, \cdot)$. From $\leq_{\text{icx}}$-monotonicity of the operator $\mathcal{T}_t$, it follows for $x \leq y$ that

$$\widetilde{\mathcal{G}}_K(t_0, x) = E^* g(\mathcal{T}_{t_{K-1}} \cdots \mathcal{T}_{t_0} x)$$
$$\leq E^* g(\mathcal{T}_{t_{K-1}} \cdots \mathcal{T}_{t_0} y) = \widetilde{\mathcal{G}}_K(t_0, y).$$

Thus, $\widetilde{\mathcal{G}}_K(t, \cdot) \in \mathcal{F}_{\text{icx}}$ and Assumption AP($g$) implies that $\mathcal{G}(t_0, \cdot) \in \mathcal{F}_{\text{icx}}$. □

Theorems 2.2 and 2.6 imply increasing convex comparison of diffusions with jumps with (13) and semimartingales $S \sim (b, c, K)$. $S^*$ has characteristics $S^* \sim (b^*, c^*, K^*)$ of the form (3).

THEOREM 2.7 (Increasing convex comparison of semimartingales to Markovian diffusions with jumps). *Let* $S^* \sim (b^*(t, S^*_{t-}), \sigma^*(t, S^*_{t-})(\sigma^*(t, S^*_{t-}))^T,$ $\lambda^{*\phi^*(t, S^*_{t-}, \cdot)})$ *be a $d$-dimensional diffusion with jumps, let* $S \sim (b, c, K)$ *be a $d$-dimensional semimartingale and assume that* $S_0 = S^*_0$. *Let* $g \in \overline{\mathcal{F}}_{\text{icx}}$ *and assume that Assumption* AP($g$) *is satisfied for* $S^*$. *Let* $b^*(t, \cdot)$ *and* $\sigma^*(t, \cdot)$ *be increasing and convex for all* $t \in [0, T]$, *as in Lemma* 2.5, *and let* $\phi^*$ *satisfy* (14) *or* (15) *($d = 1$ in the latter case). If*

$$b^i_t(\omega) \leq b^{*i}(t, S_{t-}(\omega)), \qquad c_t(\omega) \leq_{\text{psd}} (\sigma^*(t, S_{t-}(\omega)))^T \sigma^*(t, S_{t-}(\omega)),$$
$$\int_{\mathbb{R}^d} f(t, S_{t-}(\omega), x) K_{\omega, t}(dx) \leq \int_{\mathbb{R}^d} f(t, S_{t-}(\omega), x) \lambda^{*\phi^*(t, S^*_{t-}, \cdot)}(dx),$$

$\lambda \times P$-*a.e., for all* $f \in \mathbb{R}_+ \times \mathbb{R}^d \times \mathbb{R}^d \to \mathbb{R}$ *with* $f(t, s, \cdot) \in \mathcal{F}_{\text{cx}}$ *such that the integrals exist, then*

$$E g(S_T) \leq E^* g(S^*_T).$$

In the sequel, we apply the comparison results of Theorem 2.2 to Lévy processes. To that aim, we establish the propagation of order property for all order-generating function classes $\mathcal{F}$ that are considered in this paper. Also, in the following Section 3, we will establish comparison results for the finite-dimensional distributions of Lévy processes, derived by a different method.

We denote the transition probability by $P^*_{s,t}(x, B) := P^*(S^*_t - S^*_s \in B - x)$ and $P^*_{s,t}(B) := P^*_{s,t}(0, B)$.

LEMMA 2.8 [PO($g$) for spatially homogeneous processes]. *Let* $\mathcal{F}$ *be a function class as in* (1) *and let* $g \in \mathcal{F}$. *If the Markovian process* $S^*$ *has a spatially homogeneous transition function* $P^*_{t_1, t_2}$, $t_1 \leq t_2$, *then* $S^*$ *satisfies* PO($g$).

PROOF. Spatial homogeneity of $P^*_{t,T}$ implies that $\mathcal{G}(t,s) = \int g(x+s) P^*_{t,T}(dx)$. Since for $f \in \mathcal{F}$, it holds that $f_x(\cdot) := f(x + \cdot) \in \mathcal{F}$ for all $x \in \mathbb{R}^d$, the result follows from the stability under mixtures property of $\leq_\mathcal{F}$. □

Lévy processes have spatially homogeneous transition functions. Therefore, we obtain the following ordering result for Lévy processes by combining Theorem 2.2 and Lemma 2.8. In this case, the differential characteristics are deterministic and linear in time [cf. (2)].

COROLLARY 2.9 [(Directionally) convex comparison of Lévy processes]. Let $S \sim (b, \Sigma, F)_{\text{id}}, S^* \sim (b^*, \Sigma^*, F^*)_{\text{id}}$ and $S_0 = S_0^*$ be Lévy processes. Let $W$ and $W^*$ be defined as in (5).

1. Let $g \in \overline{\mathcal{F}}_{\text{idcx}}$ be bounded below or $g \in \widetilde{\mathcal{F}}_{\text{idcx}}$ and $f$ or $\nu := F \times \lambda, \nu^* := F^* \times \lambda$, assume that $(|W| * \nu)_t, (|W^*| * \nu^*)_t \in \mathcal{A}^+_{\text{loc}}$. Then the conditions

$$b \leq b^*, \qquad \Sigma^{ij} \leq \Sigma^{*ij} \qquad \forall i,j \leq d \quad \text{and} \quad F \leq_{\text{dcx}} F^*$$

imply that $Eg(S_T) \leq E^* g(S_T^*)$.

2. If $g \in \overline{\mathcal{F}}_{\text{icx}}$, then the conditions

$$b \leq b^*, \qquad \Sigma \leq_{\text{psd}} \Sigma^* \quad \text{and} \quad F \leq_{\text{cx}} F^*$$

imply that $Eg(S_T) \leq E^* g(S_T^*)$.

REMARK 2.10. According to Remark 2.3, similar statements hold true for the orders generated by $\mathcal{F}_{\text{st}}, \mathcal{F}_{\text{cx}}, \mathcal{F}_{\text{dcx}}, \mathcal{F}_{\text{icx}}, \mathcal{F}_{\text{idcx}}, \mathcal{F}_{\text{sm}}$ and $\mathcal{F}_{\text{ism}}$.

In particular, Corollary 2.9 and Remark 2.10 imply the following, interesting ordering results for multivariate normal random variables which are established in a different way in [17], Chapter 3.13:

COROLLARY 2.11 (Ordering of normal random vectors). Let $N \sim N(\mu, \Sigma)$, $N^* \sim (\mu^*, \Sigma^*)$.

1. If $\mu \leq \mu^*$ and $\Sigma = \Sigma^*$, then $N \leq_{\text{st}} N^*$.
2. If $\mu = \mu^*$ and $\Sigma \leq_{\text{psd}} \Sigma^*$, then $N \leq_{\text{cx}} N^*$.
3. If $\mu = \mu^*$ and $\Sigma^{ij} \leq \Sigma^{*ij} \forall i, j \leq d$, then $N \leq_{\text{dcx}} N^*$.
4. If $\mu \leq \mu^*$ and $\Sigma \leq_{\text{psd}} \Sigma^*$, then $N \leq_{\text{icx}} N^*$.
5. If $\mu \leq \mu^*$ and $\Sigma^{ij} \leq \Sigma^{*ij} \forall i, j \leq d$, then $N \leq_{\text{idcx}} N^*$.

**3. Comparison results for Lévy processes.** In this section, we investigate a different approach that yields comparison results for finite-dimensional distributions of Lévy processes. For this approach, we first establish an ordering result for time-homogeneous Markov processes based on a monotone separating transition kernel. We consider the finite-dimensional ordering $(S_t^1) \leq_\mathcal{F} (S_t^2)$ based on the function classes $\mathcal{F}^{(m)} := \{f := (\mathbb{R}^d)^m \to$ $\mathbb{R} : f(s_1, \ldots, s_{i-1}, \cdot, s_{i+1}, \ldots, s_m) \in \mathcal{F}, s_i \in \mathbb{R}^d, i \leq m\}$, $m, d \in \mathbb{N}$, that are componentwise in $\mathcal{F}$. A $d$-dimensional process $S^{(1)}$ is said to have smaller finite-dimensional distributions with respect to the product ordering induced by $\mathcal{F}$ than a $d$-dimensional process $S^{(2)}$, $(S_t^1) \leq_{\mathcal{F}} (S_t^2)$, if for every $m \in \mathbb{N}$ and all $0 \leq t_1 < \cdots < t_m \leq T$, it holds true that $Eg(S_{t_1}^{(1)}, \ldots, S_{t_m}^{(1)}) \leq Eg(S_{t_1}^{(2)}, \ldots, S_{t_m}^{(2)})$ for all $g \in \mathcal{F}^{(m)}$.

For time-homogeneous Markov processes, the existence of a $\leq_{\mathcal{F}}$-monotone transition kernel that separates the transition kernels of $S^{(i)}$ is sufficient to establish ordering of the finite-dimensional distributions. This is stated in the following separation result which holds true for general integral orders induced by a function class $\mathcal{F}$:

PROPOSITION 3.1 (Ordering of finite-dimensional distributions of Markov processes). *Two time-homogeneous Markov processes* $(S_t^{(1)})_{t \in [t_1, T]}$ *and* $(S_t^{(2)})_{t \in [t_1, T]}$ *with transition kernels* $Q_t^{(1)}$ *and* $Q_t^{(2)}$ *satisfy*

$$(S_t^{(1)}) \leq_{\mathcal{F}} (S_t^{(2)})$$

*if* $S_{t_1}^{(1)} \leq_{\mathcal{F}} S_{t_1}^{(2)}$ *and if a family* $(Q_t)$ *of* $\leq_{\mathcal{F}}$-*monotone transition kernels exists such that*

$$Q_t^{(1)}(x, \cdot) \leq_{\mathcal{F}} Q_t(x, \cdot) \leq_{\mathcal{F}} Q_t^{(2)}(x, \cdot) \qquad \text{for all } x \text{ and all } t \in [t_1, T].$$

PROOF. The proof uses arguments similar to those used in the proof of Theorem 5.2.15 in [17]. We consider the case of dimension $m = 2$. For $f \in \mathcal{F}^{(2)}$, $\leq_{\mathcal{F}}$-monotonicity of $Q_{t_2-t_1}(s_1, \cdot)$ is equivalent to $g(s_1) := \int f(s_1, s_2) Q_{t_2-t_1}(s_1, ds_2) \in \mathcal{F}$. As $f(s_1, \cdot) \in \mathcal{F}$, the ordering of the transition kernels implies for $f^{(i)}(s_1) := \int f(s_1, s_2) Q_{t_2-t_1}^{(i)}(s_1, ds_2)$ that $f^{(1)}(s_1) \leq g(s_1) \leq f^{(2)}(s_1)$. Hence, $S_{t_1}^{(1)} \leq_{\mathcal{F}} S_{t_1}^{(2)}$ implies that

$$Ef(S_{t_1}^{(1)}, S_{t_2}^{(1)}) = \int f^{(1)}(s_1) P^{S_{t_1}^{(1)}}(ds_1) \leq \int g(s_1) P^{S_{t_1}^{(1)}}(ds_1)$$

$$\leq \int g(s_1) P^{S_{t_1}^{(2)}}(ds_1) \leq \int f^{(2)}(s_1) P^{S_{t_1}^{(2)}}(ds_1) = Ef(S_{t_1}^{(2)}, S_{t_2}^{(2)}).$$

The result for $m > 2$ follows by induction. □

In order to apply Proposition 3.1 to Lévy processes, we consider, in the first step, the particular case of compound Poisson processes. We make use of a representation of compound Poisson processes as random sum processes with a Poisson number of summands, which allows comparison results with respect to the order-generating functions $\mathcal{F}$ from (1) by a natural coupling argument. In the second step, we extend the comparison result to pure jump

Lévy processes with infinite Lévy measure. We truncate the corresponding Lévy measures around the origin and give ordering conditions that imply ordering with respect to the order-generating functions in (1) of the induced compound Poisson processes. By an approximation argument, this implies ordering for the limit process.

For finite measures $M^{(i)}$, $i = 1, 2$, we write $M^{(1)} \leq_\mathcal{F} M^{(2)}$ to denote that $\int f(x) M^{(1)}(dx) \leq \int f(x) M^{(2)}(dx)$ for all $f \in \mathcal{F}$.

LEMMA 3.2 (Ordering of compound Poisson processes). *Let $S^{(i)} \sim (b^{(i)}(0), 0, F^{(i)})_0$, let $E|S_1^{(i)}| < \infty$, and assume that the Lévy measures $F^{(i)}$ have the same finite total mass $\lambda$. If $F^{(1)} \leq_\mathcal{F} F^{(2)}$ holds true for:*

1. $\mathcal{F} \in \{\mathcal{F}_{cx}, \mathcal{F}_{dcx}, \mathcal{F}_{sm}\}$ *and* $b^{(1)}(0) = b^{(2)}(0)$

*or for*

2. $\mathcal{F} \in \{\mathcal{F}_{st}, \mathcal{F}_{icx}, \mathcal{F}_{idcx}, \mathcal{F}_{ism}\}$ *and* $b^{(1)}(0) \leq b^{(2)}(0)$,

*then* $(S_t^{(1)}) \leq_\mathcal{F} (S_t^{(2)})$.

PROOF. Let $t \in [0, T]$.

1. Let $\mathcal{F} \in \{\mathcal{F}_{cx}, \mathcal{F}_{dcx}, \mathcal{F}_{sm}\}$. Then $S_t^{(i)}$ are representable as random sum processes

$$S_t^{(i)} = b^{(i)}(0)t + \sum_{j=1}^{N_t} X_j^{(i)}, \qquad i = 1, 2,$$

where $(X_j^{(i)})$ are i.i.d. with distribution $R^{(i)}(dx) = \frac{1}{\lambda} F^{(i)}(dx)$ and are independent of the Poisson process $N_t$ with intensity $\lambda$. Thus, $S^{(1)}$ and $S^{(2)}$ are naturally coupled by the same Poisson process $N$. As $\|F^{(i)}\| = \lambda$, condition $F^{(1)} \leq_\mathcal{F} F^{(2)}$ implies that $X_j^{(1)} \leq_\mathcal{F} X_j^{(2)}$ $\forall j$. From the stability of $\leq_\mathcal{F}$ with respect to convolutions and mixtures, the coupling of $S^{(1)}$ and $S^{(2)}$ implies, for $f \in \mathcal{F}$, that,

$$Ef\left(\sum_{j=1}^{N_t} X_j^{(1)}\right) = E^{N_t} Ef\left(\sum_{j=1}^n X_j^{(1)}\right) \leq E^{N_t} Ef\left(\sum_{j=1}^n X_j^{(2)}\right)$$
$$= E^{N_t} Ef\left(\sum_{j=1}^n X_j^{(2)}\right) = Ef\left(\sum_{j=1}^{N_t} X_j^{(2)}\right),$$

where $E^{N_t}$ denotes the expectation with respect to the distribution of $N_t$. Also, from convolution stability, it follows that $S_t^{(1)} \leq_\mathcal{F} S_t^{(2)}$.

2. Let $\mathcal{F} \in \{\mathcal{F}_{\text{st}}, \mathcal{F}_{\text{icx}}, \mathcal{F}_{\text{idcx}}, \mathcal{F}_{\text{ism}}\}$. Similarly to the first part of the proof, it follows that $\sum_{j=1}^{N_t} X_j^{(1)} \leq_{\mathcal{F}} \sum_{j=1}^{N_t} X_j^{(2)}$. As all functions in $\mathcal{F}$ are increasing and $b^{(2)}(0) - b^{(1)}(0) \geq 0$, it follows that

$$S_t^{(1)} = b^{(1)}(0)t + \sum_{j=1}^{N_t} X_j^{(1)} \leq_{\mathcal{F}} b^{(1)}(0)t + \sum_{j=1}^{N_t} X_j^{(2)}$$

$$\leq_{\mathcal{F}} b^{(2)}(0)t + \sum_{j=1}^{N_t} X_j^{(2)} = S_t^{(2)}.$$

3. It remains to prove that the finite-dimensional distributions are also ordered. By Proposition 3.1, we must establish that $Q_t^{(1)}(x, \cdot) = P(S_t^{(1)} \in \cdot \mid S_0^{(1)} = x)$ is $\leq_{\mathcal{F}}$-monotone. To that end, it suffices to prove that for $f \in \mathcal{F}$, the function

$$f_{Q_t^{(1)}}(x) := \int f(y) Q_t^{(1)}(x, dy)$$

belongs to $\mathcal{F}$. From spatial homogeneity, it follows that $f_{Q_t^{(1)}}$ is of the form $f_{Q_t^{(1)}}(x) = \int f(x+y) Q_t^{(1)}(dy)$, where $Q_t^{(1)}(dy) := Q_t^{(1)}(0, dy)$. For $\mathcal{F} = \mathcal{F}_{\text{cx}}$ and $\alpha \in (0,1)$, this implies that

$$f_{Q_t^{(1)}}(\alpha x + (1-\alpha)z) = \int f(\alpha(x+y) + (1-\alpha)(z+y)) Q_t^{(1)}(dy)$$

$$\leq \alpha f_{Q_t^{(1)}}(x) + (1-\alpha) f_{Q_t^{(1)}}(z),$$

thus $f_{Q_t^{(1)}} \in \mathcal{F}_{\text{cx}}$. Similarly, it is established for $\mathcal{F} \in \{\mathcal{F}_{\text{st}}, \mathcal{F}_{\text{icx}}, \mathcal{F}_{\text{idcx}}, \mathcal{F}_{\text{ism}}, \mathcal{F}_{\text{dcx}}, \mathcal{F}_{\text{sm}}\}$ that $f_{Q_t^{(1)}} \in \mathcal{F}$ holds true. □

We now turn to ordering results on Lévy processes with infinite Lévy measure. We truncate the corresponding Lévy measures around zero, use ordering results from Lemma 3.2 and then obtain orderings for the limit processes by weak convergence.

Let $F^{(i)}$ be Lévy measures with infinite total mass. For sequences $\underline{\varepsilon}_n^{(i)} \uparrow 0$ and $\overline{\varepsilon}_n^{(i)} \downarrow 0$, we define the truncated Lévy measures $F_n^{(i)}$ by

(19) $$F_n^{(i)}(dx) := \mathbb{1}_{(\underline{\varepsilon}_n^{(i)}, \overline{\varepsilon}_n^{(i)})^c}(x) F^{(i)}(dx).$$

For $X_n^{(i)} \sim (0, 0, F_n^{(i)})_0$, we introduce the modified truncated Lévy measures $\overline{F}_n^{(i)}$ defined as

(20) $$\overline{F}_n^{(i)}(dx) := F_n^{(i)}(dx) + (\|F_n^{(3-i)}\| - \|F_n^{(i)}\|) \mathbb{1}_{\{\|F_n^{(3-i)}\| \geq \|F_n^{(i)}\|\}} \delta_{\{0\}}(dx).$$

Then $X_n^{(i)} \sim (0, 0, \overline{F}_n^{(i)})_0$ and $\|\overline{F}_n^{(1)}\| = \|\overline{F}_n^{(2)}\|$.

THEOREM 3.3 (Convex-type comparison of Lévy processes). *Let $\mathcal{F} \in \{\mathcal{F}_{\text{cx}}, \mathcal{F}_{\text{dcx}}, \mathcal{F}_{\text{sm}}\}$. Let $S_t^{(i)}$, $i = 1, 2$, be Lévy processes with $E|S_1^{(i)}| < \infty$ and $S^{(i)} \sim (ES_1^{(i)}, 0, F^{(i)})_{\text{id}}$, where the Lévy measures $F^{(i)}$ have infinite total mass. Let $\underline{\varepsilon}_n^{(i)} \uparrow 0$ and $\overline{\varepsilon}_n^{(i)} \downarrow 0$ be sequences such that for $F_n^{(i)}$ given in* (19), *it holds true that $\int x F_n^{(1)}(dx) = \int x F_n^{(2)}(dx)$.*

*If $ES_1^{(1)} = ES_1^{(2)}$ and if the modified truncated Lévy measures $\overline{F}_n^{(i)}$ defined in* (20) *are ordered by*

$$(21) \qquad \overline{F}_n^{(1)} \leq_{\mathcal{F}} \overline{F}_n^{(2)} \qquad \forall\, n \in \mathbb{N},$$

*then $(S_t^{(1)}) \leq_{\mathcal{F}} (S_t^{(2)})$.*

PROOF. Let $S_n^{(i)} \sim (ES_1^{(i)}, 0, \overline{F}_n^{(i)})_{\text{id}} = (b^{(i)}(0), 0, \overline{F}_n^{(i)})_0$, where $b_n^{(i)}(0) := ES_1^{(i)} - \int x F_n^{(i)}(dx)$. Then $b_n^{(1)}(0) = b_n^{(2)}(0)$ and due to Lemma 3.2, the ordering on the modified truncated Lévy measures (21) implies that $(S_{n,t}^{(1)}) \leq_{\mathcal{F}} (S_{n,t}^{(2)})$. As $ES_{n,t}^{(i)} = ES_t^{(i)}$, it remains to show that $(S_{n,t}^{(i)}) \xrightarrow{\mathcal{D}} (S_t^{(i)})$ (see [17], Lemma 3.4.5, Theorems 3.4.6, 3.12.7 and 3.3.12). A sufficient condition for functional weak convergence is convergence of the differential characteristics, as specified in [13], Corollary VII.3.6. For the drift component, we have $b_n^{(i)}(h) = b^{(i)}(h)$, for all $n$. For the modified Gaussian characteristics, it follows that

$$\begin{aligned}
\widetilde{c}_n^{(i)k,\ell} &= \int h^k(x) h^\ell(x) F_n^{(i)}(dx) \\
&= \int h^k(x) h^\ell(x) \mathbb{1}_{(\underline{\varepsilon}_n^{(i)}, \overline{\varepsilon}_n^{(i)})^c}(x) F^{(i)}(dx) \\
&\longrightarrow \int h^k(x) h^\ell(x) F^{(i)}(dx) = \widetilde{c}^{(i)k,\ell},
\end{aligned}$$

due to the Lebesgue theorem. The appropriate convergence of the Lévy measures is with respect to $\widetilde{g} \in C_2(\mathbb{R}^d) := \{f : \mathbb{R}^d \to \mathbb{R},\, f$ is bounded, continuous and has value 0 around 0$\}$. It is obvious that there exists an $N \in \mathbb{N}$ such that

$$\begin{aligned}
F_n^{(i)}(g) &= \int g(x) \mathbb{1}_{(-\varepsilon_n^{(i)}, \varepsilon_n^{(i)})^c}(x) F^{(i)}(dx) \\
&= \int g(x) F^{(i)}(dx)
\end{aligned}$$

for all $n \geq N$. Therefore, Jacod and Shiryaev ([13], Corollary VII.3.6) implies functional convergence $S_n^{(i)} \xrightarrow{\mathcal{D}} S^{(i)}$, thus weak convergence of the finite-dimensional distributions, and the result follows. $\square$

From Theorem 3.3, it follows that if a Lévy process $S^{(2)}$ has higher activity in every jump height than some further Lévy process $S^{(1)}$, then $(S_t^{(1)}) \leq_{\text{cx}} (S_t^{(2)})$.

COROLLARY 3.4 (Domination criterion). *Let $S_t^{(i)}$, $S^{(i)} \sim (ES_1^{(i)}, 0, F^{(i)})_{\text{id}}$, $i = 1, 2$, be one-dimensional Lévy processes with $E|S_1^{(i)}| < \infty$ and assume that $F^{(i)}$ have densities $f^{(i)}$ such that $\int_A |x| F^{(i)}(dx) = \infty$ for $A = (-1, 0)$ and $A = (0, 1)$. If $ES_1^{(1)} = ES_1^{(2)}$ and*

$$(22) \qquad 0 < f^{(1)}(x) \leq f^{(2)}(x) \qquad \forall\, x \in \mathbb{R},$$

*then $(S_t^{(1)}) \leq_{\text{cx}} (S_t^{(2)})$.*

PROOF. Let $\underline{\varepsilon}_n^{(i)} \uparrow 0$, $i = 1, 2$, be sequences with $\underline{\varepsilon}_n^{(1)} \leq \underline{\varepsilon}_n^{(2)}$, for all $n \in \mathbb{N}$. Due to the assumption $\int_A |x| F^{(i)}(dx) = \infty$ for $A = (-1, 0)$ and $A = (0, 1)$, there are sequences $\overline{\varepsilon}_n^{(i)} \downarrow 0$, $i = 1, 2$, such that $E_n^{(i)} := \int x F_n^{(i)}(dx) = 0$. For the finite measures $\overline{F}_n^{(i)}$ given in (20), define the Lévy distribution functions $\overline{F}_n^{(i)}(x) := \overline{F}_n^{(i)}((-\infty, x])$. We establish that for all $n \in \mathbb{N}$, there is an $x_n \in \mathbb{R}$ with

$$(23) \quad \overline{F}_n^{(1)}(x) \leq \overline{F}_n^{(2)}(x) \quad \forall\, x \leq x_n \quad \text{and} \quad \overline{F}_n^{(1)}(x) \geq \overline{F}_n^{(2)}(x) \quad \forall\, x \geq x_n.$$

As $E_n^{(i)} = 0$, the cut criterion by Karlin and Novikov (cf. [17], Theorem 1.5.17) implies $\overline{F}_n^{(1)} \leq_{\text{cx}} \overline{F}_n^{(2)}$ for all $n \in \mathbb{N}$ and the result follows from Theorem 3.3.

From $f^{(1)}(x) \leq f^{(2)}(x)$ and $\underline{\varepsilon}_n^{(1)} \leq \underline{\varepsilon}_n^{(2)}$, it follows that $\overline{F}_n^{(1)}(x) \leq \overline{F}_n^{(2)}(x)$ for all $x < 0$ and it remains to prove that the cutting point $x_n$ is in the positive half-axis including the origin. The Lévy distribution functions $\overline{F}_n^{(i)}$ depend on $\overline{\varepsilon}_n^{(i)}$ and $\|F_n^{(i)}\|$. The following three cases are possible:

1. $\overline{\varepsilon}_n^{(1)} \leq \overline{\varepsilon}_n^{(2)}$, $\|F_n^{(2)}\| \geq \|F_n^{(1)}\|$,
2. $\overline{\varepsilon}_n^{(1)} \leq \overline{\varepsilon}_n^{(2)}$, $\|F_n^{(2)}\| \leq \|F_n^{(1)}\|$,
3. $\overline{\varepsilon}_n^{(1)} \geq \overline{\varepsilon}_n^{(2)}$.

Observe that $\overline{\varepsilon}_n^{(1)} \geq \overline{\varepsilon}_n^{(2)}$ implies $\|F_n^{(2)}\| \geq \|F_n^{(1)}\|$.

We establish (23) for the first case; the other cases are similar. From $\|F_n^{(2)}\| \geq \|F_n^{(1)}\|$, it follows that the modified truncated Lévy measures $\overline{F}_n^{(i)}$ are of the form $\overline{F}_n^{(2)} = F_n^{(2)}$ and $\overline{F}_n^{(1)}(dx) = F_n^{(1)}(dx) + (\|F_n^{(2)}\| - \|F_n^{(1)}\|)\delta_{\{0\}}(dx)$ and that one of the following two cases holds true:

$$\overline{F}_n^{(1)}(0) \leq \overline{F}_n^{(2)}(0) \quad \text{or} \quad \overline{F}_n^{(1)}(0) > \overline{F}_n^{(2)}(0).$$

The first case implies that $\overline{F}_n^{(1)}(x) \leq \overline{F}_n^{(2)}(x)$ for all $x \in [0, \overline{\varepsilon}_n^{(1)}]$. From $\|\overline{F}_n^{(i)}\| = \lambda_n$ and $f^{(1)}(x) \leq f^{(2)}(x)$ for all $x \in \mathbb{R}$, it follows that $\overline{F}_n^{(1)}(x) \geq \overline{F}_n^{(2)}(x)$ for all $x \geq \overline{\varepsilon}_n^{(2)}$. As $\overline{F}_n^{(2)}$ is constant on $[\overline{\varepsilon}_n^{(1)}, \overline{\varepsilon}_n^{(2)}]$ and $\overline{F}_n^{(1)}$ is strongly monotone in this interval, it follows that (23) holds true for one and only one $x_n \in [\overline{\varepsilon}_n^{(1)}, \overline{\varepsilon}_n^{(2)}]$.

In the second case, it follows from $\overline{F}_n^{(1)}(0) > \overline{F}_n^{(2)}(0)$ and $\overline{\varepsilon}_n^{(1)} \leq \overline{\varepsilon}_n^{(2)}$ that $\overline{F}_n^{(1)}(x) \geq \overline{F}_n^{(2)}(x)$ for all $x \in [0, \overline{\varepsilon}_n^{(2)}]$. The condition $f^{(1)}(x) \leq f^{(2)}(x)$ for all $x \in \mathbb{R}$, together with $\|\overline{F}_n^{(i)}\| = \lambda_n$, implies that $\overline{F}_n^{(1)}(x) \geq \overline{F}_n^{(2)}(x)$ for $x > \overline{\varepsilon}_n^{(2)}$. Therefore, (23) holds true for $x_n = 0$. □

Corollary 3.4 can be applied to compare one-dimensional normal inverse Gaussian processes in the shape and scaling parameters $\alpha$ and $\delta$. The Lévy density of an NIG=NIG$(\alpha, \beta, \delta, \mu)$-distributed random variable $S$ is given by

$$(24) \qquad f_{\alpha,\beta,\delta}(x) = \frac{\delta \alpha K_1(\alpha|x|)e^{\beta x}}{\pi |x|},$$

where $K_1$ denotes the modified Bessel function of the third kind with index 1, $\alpha > 0$, $0 \leq |\beta| \leq \alpha$, and $\delta > 0$ and $S$ has expectation $ES = \mu + \frac{\delta \beta}{\sqrt{\alpha^2 - \beta^2}}$. In the following example, where we give an ordering in the parameters $\alpha$ and $\delta$, respectively, we need equality of the expectations $ES^{(i)}$, which is obtained by choosing $\mu^{(i)}$ appropriately.

EXAMPLE 3.5 (Convex ordering of one-dimensional NIG processes).

1. Ordering in $\alpha$. Let $S^{(i)}, i = 1, 2$, be NIG processes with $S^{(i)} \sim (ES_1^{(i)}, 0, f^{(i)}(x)\,dx)$, where for $f_{\alpha,\beta,\delta}$ given in (24), we define $f^{(i)}(x) := f_{\alpha^{(i)},\beta,\delta}(x)$. If $\alpha^{(1)} \leq \alpha^{(2)}$ and $ES_1^{(1)} = ES_1^{(2)}$, then $(S_t^{(2)}) \leq_{\mathrm{cx}} (S_t^{(1)})$.
2. Ordering in $\delta$. Let $S^{(i)}, i = 1, 2$, be NIG processes with $S^{(i)} \sim (ES_1^{(i)}, 0, f^{(i)}(x)\,dx)$, where for $f_{\alpha,\beta,\delta}$ given in (24), we define $f^{(i)}(x) := f_{\alpha,\beta,\delta^{(i)}}(x)$. If $\delta^{(1)} \leq \delta^{(2)}$ and $ES_1^{(1)} = ES_1^{(2)}$, then $(S_t^{(1)}) \leq_{\mathrm{cx}} (S_t^{(2)})$.

PROOF. 1. For fixed $x > 0$, consider $g(\alpha) := f_{\alpha,\beta,\delta}(x)$. Then $g'(\alpha) = -\frac{\delta e^{\beta x}}{\pi} \alpha \cdot K_0(\alpha x) \leq 0$, that is, $f^{(1)}(x) \geq f^{(2)}(x)$ for all $x \in \mathbb{R}_+$. For fixed $x < 0$, we similarly obtain $g'(\alpha) = -\frac{\delta e^{\beta x}}{\pi} \alpha K_0(-\alpha x) \leq 0$, thus $f^{(1)}(x) \geq f^{(2)}(x)$, for all $x \in \mathbb{R}_-$. Corollary 3.4 implies that $(S_t^{(2)}) \leq_{\mathrm{cx}} (S_t^{(1)})$.

2. It is obvious that $\delta^{(1)} \leq \delta^{(2)}$ implies $f^{(1)}(x) \leq f^{(2)}(x)$, for all $x \in \mathbb{R}$ and, therefore, by Corollary 3.4, it follows that $(S_t^{(1)}) \leq_{\mathrm{cx}} (S_t^{(2)})$. □

The following result for increasing- and increasing convex-type orders is similar to Theorem 3.3 and Corollary 3.4:

THEOREM 3.6 (Increasing convex-type comparison of Lévy processes). *Let $\mathcal{F} \in \{\mathcal{F}_{\text{st}}, \mathcal{F}_{\text{icx}}, \mathcal{F}_{\text{idcx}}, \mathcal{F}_{\text{ism}}\}$ and assume that $S_t^{(i)}, i = 1, 2$, are Lévy processes with $E|S_1^{(i)}| < \infty$ and $S^{(i)} \sim (ES_1^{(i)}, 0, F^{(i)})_{\text{id}}$, where $\|F^{(i)}\| = \infty$. Assume that $ES_1^{(1)} \leq ES_1^{(2)}$ and let $\underline{\varepsilon}_n^{(i)} \uparrow 0$ and $\overline{\varepsilon}_n^{(i)} \downarrow 0$ be sequences such that for $F_n^{(i)}$ given in (19), it holds true that $0 \leq \int xF_n^{(2)}(dx) - \int xF_n^{(1)}(dx) \leq ES_1^{(2)} - ES_1^{(1)}$.*

*If the modified truncated Lévy measures $\overline{F}_n^{(i)}$ defined in (20) are ordered as*

$$\overline{F}_n^{(1)} \leq_\mathcal{F} \overline{F}_n^{(2)} \qquad \forall n \in \mathbb{N},$$

*then $(S_t^{(1)}) \leq_\mathcal{F} (S_t^{(2)})$.*

PROOF. For $b_n^{(i)}(0) := ES_1^{(i)} - \int xF_n^{(i)}(dx)$ and $S_n^{(i)} \sim (b_n^{(i)}(0), 0, \overline{F}_n^{(i)})_0$, it follows from Lemma 3.2 that $(S_{n,t}^{(1)}) \leq_\mathcal{F} (S_{n,t}^{(2)})$. As $ES_{n,1}^{(i)} = b_n^{(i)}(\text{id}) = ES_1^{(i)}$, it remains to prove functional weak convergence. This is similar to the proof of Theorem 3.3. □

If the processes $S^{(i)}$ have paths of finite variation, then there is a similar comparison result that does not require the first moment condition $0 \leq \int xF_n^{(2)}(dx) - \int xF_n^{(1)}(dx) \leq ES_1^{(2)} - ES_1^{(1)}$. In this case, we postulate an ordering condition on the drift terms $b^{(1)}(0) \leq b^{(2)}(0)$ of $S^{(i)}$. The $b^{(i)}(0)$ exist, due to the fact that Lévy measures of processes with paths of finite variation integrate $|x|$ around the origin.

THEOREM 3.7 (Increasing convex-type comparison of Lévy processes with paths of finite variation). *Let $\mathcal{F} \in \{\mathcal{F}_{\text{st}}, \mathcal{F}_{\text{icx}}, \mathcal{F}_{\text{idcx}}, \mathcal{F}_{\text{ism}}\}$ and assume that $S_t^{(i)}, i = 1, 2$, are Lévy processes with paths of finite variation. Assume that $E|S_1^{(i)}| < \infty$ and $S^{(i)} \sim (b^{(i)}(0), 0, F^{(i)})_0$, where $\|F^{(i)}\| = \infty$. Let $\underline{\varepsilon}_n^{(i)} \uparrow 0$ and $\overline{\varepsilon}_n^{(i)} \downarrow 0$ be sequences such that for $F_n^{(i)}$ given in (19), it holds true that $\int xF_n^{(1)}(dx) \leq \int xF_n^{(2)}(dx)$.*

*If $b^{(1)}(0) \leq b^{(2)}(0)$ and the modified truncated Lévy measures $\overline{F}_n^{(i)}$ defined in (20) are ordered as*

$$\overline{F}_n^{(1)} \leq_\mathcal{F} \overline{F}_n^{(2)} \qquad \forall n \in \mathbb{N},$$

*then $(S_t^{(1)}) \leq_\mathcal{F} (S_t^{(2)})$.*

PROOF. For $S_n^{(i)} \sim (b^{(i)}(0), 0, \overline{F}_n^{(i)})_0$, it follows from Lemma 3.2 that $(S_{n,t}^{(1)}) \leq_\mathcal{F} (S_{n,t}^{(2)})$. As $ES_{n,1}^{(i)} = b_n^{(i)}(\text{id}) = b^{(i)}(0) + \int xF_n^{(i)}(dx) \to b^{(i)}(\text{id}) = ES_1^{(i)}$, it remains to prove functional weak convergence. This is similar to the proof of Theorem 3.3. □

For the one-dimensional Theorem 3.7, a sufficient condition for the ordering with respect to $\mathcal{F}_{\text{st}}$ in terms of Lévy densities. If the Lévy density of one of the processes is dominated by the Lévy density of the other on the negative half-axis and the domination is reversed on the positive half-axis, then the first process is smaller than the second one with respect to the usual stochastic ordering.

COROLLARY 3.8. *Let $S_t^{(i)}, i = 1, 2$, be one-dimensional Lévy processes with finite variation, $E|S_1^{(i)}| < \infty$ and $S^{(i)} \sim (b^{(i)}(0), 0, F^{(i)})_0$. Assume that $F^{(i)}$ have densities $f^{(i)}$ that are monotonically increasing to infinity as $x$ tends to zero.*
*If $b^{(1)}(0) \leq b^{(2)}(0)$,*

$$f^{(1)}(x) \geq f^{(2)}(x) \quad \forall\, x \in \mathbb{R}_- \quad \text{and} \quad f^{(1)}(x) \leq f^{(2)}(x) \quad \forall\, x \in \mathbb{R}_+,$$

*then $(S_t^{(1)}) \leq_{\text{st}} (S_t^{(2)})$.*

PROOF. Let $\varepsilon_n \downarrow 0$ and let the sequences $\underline{\varepsilon}_n^{(i)} \uparrow 0$, $\overline{\varepsilon}_n^{(i)} \downarrow 0$ in (19) be given by $\overline{\varepsilon}_n^{(i)} := \varepsilon_n$, $\varepsilon_n := -\underline{\varepsilon}_n^{(i)}$, $i = 1, 2$. Then for $E_n^{(i)} := \int x F_n^{(i)}(dx)$, it follows that $E_n^{(1)} \leq E_n^{(2)}$, due to the pointwise ordering condition on the Lévy densities $f^{(i)}$. Similarly to Corollary 3.4, we establish that

$$(25) \qquad \overline{F}_n^{(1)}(x) \geq \overline{F}_n^{(2)}(x) \qquad \forall\, x \in \mathbb{R}.$$

This implies that $\overline{F}_n^{(1)} \leq_{\text{st}} \overline{F}_n^{(2)}$, for all $n \in \mathbb{N}$ and the result follows from Theorem 3.7.

From $f^{(1)}(x) \geq f^{(2)}(x)$, for all $x < 0$ and $\underline{\varepsilon}_n^{(i)} = -\varepsilon_n$ it follows that $F_n^{(1)}(x) \geq F_n^{(2)}(x)$ for all $x < 0$ and it only remains to prove the domination criterion (25) on the positive half-axis including the origin. The modified truncated Lévy measures depend on the sign of $\|F_n^{(2)}\| - \|F_n^{(1)}\|$. If $\|F_n^{(2)}\| \geq \|F_n^{(1)}\|$, then $\overline{F}_n^{(2)} = F_n^{(2)}$ and $\overline{F}_n^{(1)}(dx) = F_n^{(1)}(dx) + (\|F_n^{(2)}\| - \|F_n^{(1)}\|)\delta_{\{0\}}(dx)$. This implies that $\overline{F}^{(1)}(x) \geq \overline{F}^{(2)}(x)$, $x \in [0, \varepsilon_n]$. From $f^{(1)}(x) \leq f^{(2)}(x)$, for all $x \in \mathbb{R}_+$, it then follows that $\overline{F}^{(1)}(x) \geq \overline{F}^{(2)}(x)$ also for $x > \varepsilon_n$ and (25) follows. If $\|F_n^{(2)}\| \leq \|F_n^{(1)}\|$ then a similar consideration implies (25). □

Due to the assumption of finite variation, the previous comparison criterion is not applicable to two NIG processes with different skewness parameters $\beta^{(1)} \leq \beta^{(2)}$, although the domination criterion on the Lévy densities of Corollary 3.8 is satisfied.

**4. Ordering result for mixing-type distributions.** For some cases of interest, it is possible to obtain comparison results by using mixing-type representations. We apply this approach to the class of generalized hyperbolic (GH) distributions. Further mixing-type representations of particular interest for financial mathematical models hold for multivariate $t$-distributions and elliptically contoured distributions (see [4]). GH distributions are variance mixtures of multivariate normal distributions which have a generalized inverse Gaussian distribution as mixing distribution. For $\mu^{(i)}, \beta^{(i)} \in \mathbb{R}^d$, $\Delta^{(i)} \in M(d, \mathbb{R})$ with $\det(\Delta^{(i)}) = 1$ for $i = 1, 2$ and $N^{(i)} \sim \mathcal{N}(0, \Delta^{(i)})$, we consider the $d$-dimensional random variable

$$(26) \qquad S^{(i)} = \mu^{(i)} + X^{(i)} \Delta^{(i)} \beta^{(i)} + \sqrt{X^{(i)}} N^{(i)},$$

where $X^{(i)}$ are generalized inverse Gaussian random variables with densities

$$(27) \qquad d_{\mathrm{GIG}(\lambda,\delta,\gamma)}(x) := \left(\frac{\gamma}{\delta}\right)^\lambda \frac{1}{2K_\lambda(\delta\gamma)} x^{\lambda-1} e^{(-1/2)(\delta^2/x + \gamma^2 x)} \mathbb{1}_{\mathbb{R}_+}(x),$$

where $\delta \geq 0$, $\alpha^2 > \beta \Delta \beta^T$ and $\gamma = \sqrt{\alpha^2 - \beta \Delta \beta^T}$. Then $S$ is generalized hyperbolic distributed with parameters $d, \lambda, \alpha, \beta, \delta, \mu$ and covariance matrix $\Delta$ and we write $\mathrm{GH}(d, \lambda, \alpha, \beta, \delta, \mu)$ (cf. [1]).

The following lemma states a comparison result for GIG distributions with respect to the likelihood ratio order $\leq_{\mathrm{lr}}$, if the parameters $\lambda$, $\delta$ and $\gamma$ are ordered:

LEMMA 4.1 (Likelihood ratio ordering of GIG random variables). *Let $X^{(i)}$ be GIG distributed with density $d_{\mathrm{GIG}(\lambda^{(i)}, \delta^{(i)}, \gamma^{(i)})}(x)$.*
*If $\lambda^{(1)} \leq \lambda^{(2)}$, $\delta^{(1)} \leq \delta^{(2)}$ and $\gamma^{(1)} \geq \gamma^{(2)}$, then $X^{(1)} \leq_{\mathrm{lr}} X^{(2)}$.*

PROOF. We consider the likelihood ratio

$$\begin{aligned} g(x) &:= \frac{d_{\mathrm{GIG}(\lambda^{(1)}, \delta^{(1)}, \gamma^{(1)})}(x)}{d_{\mathrm{GIG}(\lambda^{(2)}, \delta^{(2)}, \gamma^{(2)})}(x)} \\ &= K x^{\lambda^{(1)} - \lambda^{(2)}} e^{(1/2)(((\delta^{(2)^2} - \delta^{(1)^2})/x) + (\gamma^{(2)^2} - \gamma^{(1)^2})x)} \mathbb{1}_{\mathbb{R}_+}(x), \end{aligned}$$

with $K > 0$. The first derivative of $g$ is of the form $g'(x) = K_1(x)(\lambda^{(1)} - \lambda^{(2)}) + K_2(x)(\frac{\delta^{(1)^2} - \delta^{(2)^2}}{x^2} + (\gamma^{(2)^2} - \gamma^{(1)^2}))$, where $K_1(x), K_2(x) \geq 0$. Thus, it follows from the orderings on the parameters that $g'(x) \leq 0$. □

In the following theorem we give an increasing convex comparison result of multivariate GH distributions with mixing-type representation (26). We consider the following three cases for $\beta^{(i)}$ and $\Delta^{(i)}$:

$$(28) \qquad 0 \leq \beta^{(1)} \leq \beta^{(2)}, \qquad \Delta^{(i)} = I,$$

$$\beta^{(i)} = 0, \qquad \Delta^{(1)} \leq_{\text{psd}} \Delta^{(2)}, \tag{29}$$

$$0 \leq \beta^{(1)} \leq \beta^{(2)}, \qquad \Delta^{(1)} \leq_{\text{psd}} \Delta^{(2)}, \tag{30}$$

$$0 \leq \Delta_{ij}^{(1)} \leq \Delta_{ij}^{(2)} \qquad \forall i,j \leq d.$$

THEOREM 4.2 (Increasing convex comparison of GH distributions). *Let $S^{(i)}$ be $\text{GH}(d, \lambda^{(i)}, \alpha^{(i)}, \beta^{(i)}, \delta^{(i)}, \mu^{(i)})$ distributed. If*

$$\lambda^{(1)} \leq \lambda^{(2)}, \qquad \delta^{(1)} \leq \delta^{(2)}, \qquad \alpha^{(1)} \geq \alpha^{(2)}, \qquad \mu^{(1)} \leq \mu^{(2)}$$

*and one of the cases* (28)–(30) *holds true for $\beta^{(i)}$ and $\Delta^{(i)}$, then $S^{(1)} \leq_{\text{icx}} S^{(2)}$.*

PROOF. First we prove that the conditions on the parameters of $X^{(i)}$ imply that $X^{(1)} \leq_{\text{st}} X^{(2)}$. For $\gamma^{(i)} = \sqrt{\alpha^{(i)^2} - \beta^{(i)} \Delta^{(i)} (\beta^{(i)})^T}$, it follows from the conditions on $\alpha^{(i)}$, $\Delta^{(i)}$ and $\beta^{(i)}$ in all three cases (28)–(30) that $\gamma^{(1)} \geq \gamma^{(2)}$. Due to Lemma 4.1, the ordering conditions on $\lambda^{(i)}$ and $\delta^{(i)}$ imply that $X^{(1)} \leq_{\text{lr}} X^{(2)}$ and this yields $X^{(1)} \leq_{\text{st}} X^{(2)}$.

1. Let $\beta^{(i)}$ and $\Delta^{(i)}$ satisfy (28). For $f \in \mathcal{F}_{\text{icx}}$, we define $g(x) := Ef(\mu^{(1)} + \beta^{(1)}x + \sqrt{x}N)$, with $N \stackrel{d}{=} N^{(i)}$, $i = 1, 2$. Let $0 \leq x^{(1)} \leq x^{(2)}$. As $\beta^{(1)}x + \sqrt{x}N \sim \mathcal{N}(\beta^{(1)}x, xI)$, $\beta^{(1)}x^{(1)} \leq \beta^{(1)}x^{(2)}$ and $x^{(1)}I \leq_{\text{psd}} x^{(2)}I$, it follows from Corollary 2.11 that $g$ is increasing. Then $X^{(1)} \leq_{\text{st}} X^{(2)}$ implies that

$$Ef(\mu^{(1)} + \beta^{(1)} X^{(1)} + \sqrt{X^{(1)}}N) = Eg(X^{(1)}) \leq Eg(X^{(2)})$$
$$= Ef(\mu^{(1)} + \beta^{(1)} X^{(2)} + \sqrt{X^{(2)}}N).$$

As $f$ is increasing and $X^{(2)}$ is nonnegative, it follows from $\mu^{(1)} \leq \mu^{(2)}$ and $\beta^{(1)} \leq \beta^{(2)}$ that

$$Ef(\mu^{(1)} + \beta^{(1)} X^{(2)} + \sqrt{X^{(2)}}N)$$
$$= E^{X^{(2)}} Ef(\mu^{(1)} + \beta^{(1)} x^{(2)} + \sqrt{x^{(2)}}N)$$
$$\leq E^{X^{(2)}} Ef(\mu^{(2)} + \beta^{(2)} x^{(2)} + \sqrt{x^{(2)}}N) = Ef(S^{(2)}).$$

2. Let $\beta^{(i)}$ and $\Delta^{(i)}$ satisfy (29) and assume that $f \in \mathcal{F}_{\text{icx}}$ and $x^{(1)} \geq 0$. From $\Delta^{(1)} \leq_{\text{psd}} \Delta^{(2)}$, it follows that $x^{(1)} \Delta^{(1)} \leq_{\text{psd}} x^{(1)} \Delta^{(2)}$ and, therefore, Corollary 2.11.4 implies $Ef(S^{(1)}) \leq Ef(\mu^{(2)} + \sqrt{X^{(1)}}N^{(2)})$. For $z \geq 0$, let $g(z) := Ef(\mu^{(2)} + zN^{(2)})$. As $\Delta^{(2)}$ is positive semidefinite, it follows for $z^{(1)} \leq z^{(2)}$ that $z^{(1)} \Delta^{(2)} \leq_{\text{psd}} z^{(2)} \Delta^{(2)}$ and Corollary 2.11.4 implies that $g$ is increasing. From $X^{(1)} \leq_{\text{st}} X^{(2)}$, it follows that

$$Ef(\mu^{(2)} + \sqrt{X^{(1)}}N^{(2)}) = Eg(\sqrt{X^{(1)}}) \leq Eg(\sqrt{X^{(2)}})$$
$$= Ef(\mu^{(2)} + \sqrt{X^{(2)}}N^{(2)}).$$

3. Let $\beta^{(i)}$ and $\Delta^{(i)}$ satisfy (30). The condition $\Delta^{(i)}_{ij} \geq 0$, $i,j \leq d$, implies that $\beta^{(1)} \Delta^{(i)} \beta^{(1)T} \leq \beta^{(2)} \Delta^{(i)} \beta^{(2)T}$ for $0 \leq \beta^{(1)} \leq \beta^{(2)}$. The result follows similarly to the previous ones. □

In the case $\lambda = -\frac{1}{2}$, $S^{(i)}$ is normally inverse Gaussian distributed. NIG distributed random variables are stable under convolutions:
$$\mathrm{NIG}(d, \alpha, \beta, \delta, \mu, \Delta; t) = \mathrm{NIG}(d, \alpha, \beta, t\delta, t\mu, \Delta; 1).$$
Therefore, Theorem 4.2 also implies increasing convex comparison of the of NIG processes with mixing-type representation

$$S^{(i)}_t := \mu^{(i)} t + X^{(i)}_t \Delta^{(i)} \beta^{(i)} + \sqrt{X^{(i)}_t} N^{(i)}, \tag{31}$$

where $X^{(i)}_t \sim \mathrm{GIG}(-\frac{1}{2}, t\delta^{(i)}, \gamma^{(i)})$ for any time $t > 0$.

COROLLARY 4.3 (Increasing convex comparison of NIG processes). *Let $S^{(i)}_t$ be $\mathrm{NIG}(d, \alpha^{(i)}, \beta^{(i)}, \delta^{(i)}, \mu^{(i)}, \Delta^{(i)}; t)$ processes. If*
$$\delta^{(1)} \leq \delta^{(2)}, \qquad \alpha^{(1)} \geq \alpha^{(2)}, \qquad \mu^{(1)} \leq \mu^{(2)}$$
*and one of the cases* (28)–(30) *holds true for $\beta^{(i)}$ and $\Delta^{(i)}$, then $(S^{(1)}_t) \leq_{\mathrm{icx}} (S^{(2)}_t)$.*

PROOF. As in Theorem 4.2, it follows from the ordering conditions that $S^{(1)}_t \leq_{\mathrm{icx}} S^{(2)}_t$ for all $t \in [0, T]$. The result on finite-dimensional distributions follows from Proposition 3.1. □

## APPENDIX

PROOF OF LEMMA 2.1. Let $t \in (0, T]$ and assume that

$$\int_{[0,t] \times \mathbb{R}^d} |\Lambda \mathcal{G}(u, S^*_{u-}(\omega^*), y)| \mu^*(\omega^*; du, dy) \in \mathcal{A}^+_{\mathrm{loc}}. \tag{32}$$

Itô's lemma implies that $\mathcal{G}(t, S^*_t)$ is a semimartingale with evolution

$$\mathcal{G}(t, S^*_t) = \mathcal{G}(0, S^*_0) + M_t + \int_{[0,t] \times \mathbb{R}^d} (\Lambda \mathcal{G})(u, S^*_{u-}, y) \mu^*(du, dy)$$
$$+ \int_{[0,t]} \left\{ \mathrm{D}_t \mathcal{G}(u, S^*_{u-}) + \sum_{i \leq d} \mathrm{D}_i \mathcal{G}(u, S^*_{u-}) b^{*i}(u, S^*_{u-}) \right.$$
$$\left. + \tfrac{1}{2} \sum_{i,j \leq d} \mathrm{D}^2_{ij} \mathcal{G}(u, S^*_{u-}) c^{*ij}(u, S^*_{u-}) \right\} du,$$

where $M_t := \sum_{i \leq d} \int_{[0,t]} \mathrm{D}_i \mathcal{G}(u, S^*_{u-}) \, dM^{*i}_u$ is a one-dimensional local $(\mathcal{A}^*_t)$-martingale under $P^*$ and $M^{*i}_t$ denotes the $i$th component of the martingale part of $S^*$. As (32) holds, there is a local martingale $\overline{M}_t$ such that

$$\begin{aligned}\mathcal{G}(t, S^*_t) &= \mathcal{G}(0, S^*_0) + M_t + \overline{M}_t \\ &\quad + \int_{[0,t]} \Bigg\{ \mathrm{D}_t \mathcal{G}(u, S^*_{u-}) + \sum_{i \leq d} \mathrm{D}_i \mathcal{G}(u, S^*_{u-}) b^{*i}(u, S^*_{u-}) \\ &\qquad\qquad + \tfrac{1}{2} \sum_{i,j \leq d} \mathrm{D}^2_{ij} \mathcal{G}(u, S^*_{u-}) c^{*ij}(u, S^*_{u-}) \\ &\qquad\qquad\qquad + \int_{\mathbb{R}^d} (\Lambda \mathcal{G})(u, S^*_{u-}, y) K^*_u(S^*_{u-}, dy) \Bigg\} du \\ &=: \mathcal{G}(0, S^*_0) + M_t + \overline{M}_t + V_t,\end{aligned}$$

where we denote the Lebesgue integral by $V_t$. As $\mathcal{G}(t, S^*_t), M_t$ and $\overline{M}_t$ are local $(\mathcal{A}^*_t)$-martingales under $P^*$, $V_t$ is a local martingale starting at zero, and as $V_t$ is of finite variation, it follows that $V \equiv 0$. Therefore, we obtain

$$\mathrm{D}_t \mathcal{G}(t, S^*_{t-}) + \sum_{i \leq d} \mathrm{D}_i \mathcal{G}(t, S^*_{t-}) b^{*i}(t, S^*_{t-})$$
$$+ \tfrac{1}{2} \sum_{i,j \leq d} \mathrm{D}^2_{ij} \mathcal{G}(t, S^*_{t-}) c^{*ij}(t, S^*_{t-}) + \int_{\mathbb{R}^d} (\Lambda \mathcal{G})(t, S^*_{t-}, y) K^*_t(S^*_{t-}, dy) = 0$$

and (4) follows.

It remains to prove that $\mathcal{G}(t, \cdot) \in \mathcal{F}_{\mathrm{cx}}$ implies the integrability condition (32). As $\mathcal{G}(t, S^*_t)$ is a local martingale under $Q^*$, it is a special semimartingale. The process

$$\int_{[0,t]} \mathrm{D}_t \mathcal{G}(u, S^*_{u-}) \, du + \int_{[0,t]} \sum_{i \leq d} \mathrm{D}_i \mathcal{G}(u, S^*_{u-}) \, dA^{*i}$$
$$+ \tfrac{1}{2} \sum_{i,j \leq d} \int_{[0,t]} \mathrm{D}^2_{ij} \mathcal{G}(u, S^*_{u-}) \, dC^{*ij}_u$$

is predictable and of finite variation and is therefore in $\mathcal{A}_{\mathrm{loc}}$. Due to a representation result for special semimartingales (see [13], Proposition I.4.23), it follows that $(W^* * \mu^*)_t \in \mathcal{A}_{\mathrm{loc}}$. As convexity of $\mathcal{G}(t, \cdot)$ implies $\Lambda \mathcal{G}(t, S^*_{t-}, y) \geq 0$, (32) follows. $\square$

**Acknowledgment.** The authors thank a reviewer for the careful reading of the manuscript and for several valuable comments.

Department of Mathematical Stochastics
Eckerstr. 1
D-79104 Freiburg
Germany
E-mail: bergen@stochastik.uni-freiburg.de
ruschen@stochastik.uni-freiburg.de